\documentclass[twocolumn,nameyear]{autart}

\usepackage{amsmath,amssymb}
\usepackage{mathtools}
\usepackage{mathrsfs}
\usepackage{graphicx}
\usepackage{wrapfig}
\usepackage{xcolor}
\newtagform{revision}{\color{blue}(}{)}
\makeatletter
\newcommand*{\wrapresultname}{%
  \def\@opargbegintheorem##1##2##3{%
    \@begintheorem{##1}{##2}%
    {\normalfont\bfseries(##3)}\enspace}}
\makeatother
\usepackage[authoryear,round]{natbib}
\usepackage[hidelinks]{hyperref}
\renewcommand*{\theHsection}{\arabic{section}}

\begin{document}

\begin{frontmatter}
\runtitle{Extremum Seeking for Static Maps via Boundary Control of an Euler-–Bernoulli PDE}
\title{Extremum Seeking for Static Maps via Boundary Control \\of an Euler--Bernoulli PDE}

\author[ufsc]{Paulo Henrique F. Biazetto}\ead{paulo.biazetto@posgrad.ufsc.br},
\author[ufsc]{Gustavo A. de Andrade}\ead{gustavo.artur@ufsc.br},\break
\author[uerj]{Tiago Roux Oliveira}\ead{tiagoroux@uerj.br},
\author[ucsd]{Miroslav Krstic}\ead{krstic@ucsd.edu}

\address[ufsc]{Department of Automation and Systems Engineering, Federal University of Santa Catarina, 88040-370 Florian\'opolis, Brazil}
\address[uerj]{Department of Electronics and Telecommunication Engineering, State University of Rio de Janeiro, Rio de Janeiro 20550-900, Brazil}
\address[ucsd]{Department of Mechanical and Aerospace Engineering, University of California, San Diego, CA 92093-0411, USA}

\begin{abstract}
This paper develops an extremum seeking (ES) controller for an unknown locally quadratic scalar map with actuation dynamics governed by the one-dimensional Euler–Bernoulli beam PDE. The displacement at the sliding end of the beam serves as the input to the map, while position and moment inputs are applied at the opposite end. The controller combines backstepping compensation of the known PDE dynamics with online gradient and Hessian estimates. The averaged closed-loop system is exponentially stable under an explicit filter bandwidth condition independent of the unknown Hessian. For a fixed nonzero probing amplitude $a$ within the quadratic region and sufficiently high frequency $\omega$, the error relative to the prescribed probing motion converges locally and exponentially to a periodic trajectory of size $\mathcal O(\omega^{-1/2})$. The sliding-end displacement therefore oscillates about the optimizer with an ultimate error bound of $|a|+\mathcal O(\omega^{-1/2})$. The position and moment probes required at the actuated end grow as $e^{\sqrt{\omega}}$ and $\omega e^{\sqrt{\omega}}$, respectively. Numerical simulations demonstrate the effectiveness of the proposed ES controller.
\end{abstract}

\begin{keyword}
Boundary control \sep Euler--Bernoulli beam \sep extremum seeking \sep optimization \sep partial differential equations
\end{keyword}
\end{frontmatter}

\section{Introduction}

Boundary control of the Euler--Bernoulli (EB) beam equation has been widely studied in the control of vibrations in flexible structures~\citep{Guo2002}. The literature includes output-feedback stabilization in the presence of boundary input disturbances~\citep{Jin2015}, exponential stabilization using boundary springs and dampers~\citep{Baysal2024}, and disturbance rejection under input saturation~\citep{Zhong2024}. Other contributions address adaptive event-triggered boundary control under uncertainties and constraints~\citep{Sun2025} and asymptotic stabilization with saturated boundary feedback~\citep{Ayadi2026}. The configuration considered here has a sliding boundary at one end and position and moment actuation at the other. For this configuration, \citet{Smyshlyaev2009} exploited the fact that a complex combination of velocity and curvature satisfies a Schr\"odinger PDE to design a backstepping controller with a prescribed exponential decay rate. In the present design, this representation yields a Schr\"odinger PDE coupled to the ODEs governing the boundary errors. The resulting cascade is related to the class studied by~\citet{Ren2013}, who used backstepping to stabilize linear ODEs with Schr\"odinger input dynamics.

Beyond regulation to a prescribed configuration, the desired operating point may be specified by the extremum of an unknown performance map. Extremum seeking (ES) addresses this problem through online measurements and probing signals, without requiring an analytical expression for the map~\citep{Krsti2000StabilityOE}. When the actuator and the map input are separated by distributed dynamics, the control design must account for the propagation of both the feedback action and the probing signals. The compensation of arbitrarily long input delays in~\citet{Oliveira2017ExtremumSF}, interpreted through transport PDEs, provided a foundation for extensions to other infinite-dimensional systems. These developments include PDE--PDE cascades~\citep{cascades}, with systematic treatments of ES through delays and PDEs available in~\citet{book2022,OliveiraKrstic2023}.

Related contributions by the authors have addressed ES for wave PDEs with Kelvin--Voigt damping~\citep{Silva2024KelvinVoigt} and for complex-valued Schr\"odinger PDEs~\citep{Biazetto2026Schrodinger}. The former combines boundary backstepping with gradient and Hessian estimation to compensate the wave dynamics. The latter uses a real-coordinate representation of the complex boundary state and a Newton-based ES law with inverse-Hessian estimation. In the present problem, the optimization variable is the real displacement at the uncontrolled end, and the control design must also account for the spatially constant displacement component and the dynamic realization of the two real boundary inputs.

This paper proposes a boundary ES controller that uses position and moment actuation at one end of an Euler--Bernoulli beam to steer the displacement at the opposite, sliding end toward the optimizer of an unknown, locally quadratic scalar map. Using the known PDE dynamics, we construct boundary probing signals that induce a prescribed sinusoidal perturbation at the map input. We then demodulate the measured map output to obtain online estimates of the gradient and Hessian. These estimates enter a backstepping feedback law derived by expressing the distributed error dynamics in Schr\"odinger form and retaining the sliding-end displacement error as a separate scalar variable.

The main contributions are as follows.
\begin{itemize}
\item We apply the boundary feedback of~\citet{Smyshlyaev2009} directly and filter only the correction formed from the gradient and Hessian estimates. This placement preserves the prescribed decay rate of the homogeneous PDE modes and yields a closed finite-dimensional subsystem for the filter and the boundary-error coordinates. With the filter dynamics retained, we prove exponential stability of the averaged closed loop whenever the filter bandwidth satisfies $\overline c>\mu_p/\mu_d$. The positive gains $\mu_p$ and $\mu_d$ are chosen by design, so this condition is independent of the unknown Hessian.

\item We derive a high-frequency bound on the response of the averaged system to the oscillatory ES terms. The bound holds at the PDE modal frequencies despite the nonanalyticity of the closed-loop semigroup. For a fixed nonzero probing amplitude $a$ within the quadratic region and sufficiently large frequency $\omega$, we prove local exponential convergence of the full closed-loop error, measured relative to the prescribed probing motion, to a periodic solution of size $\mathcal O(\omega^{-1/2})$. The sliding-end displacement consequently has an ultimate error bound of $|a|+\mathcal O(\omega^{-1/2})$. The explicit amplitudes of the position and moment probes quantify the boundary actuation required as $\omega$ increases.
\end{itemize}

The preliminary design in~\citet{Biazetto2025EulerBernoulliCDC} introduced the boundary probing signals and the two-step backstepping construction for this configuration. In that design, a low-pass filter acted on each complete boundary feedback law. Here only the map-dependent terms are filtered. Numerical simulations illustrate the effectiveness of the proposed controller.

The remainder of the paper is organized as follows. Section~\ref{section:problem} introduces the PDE model and the optimization objective. Section~\ref{section:extremum_seeking_design} presents the probing signals, error dynamics, and boundary controller. Stability and convergence are analyzed in Section~\ref{section:stability}. Section~\ref{section:simulation} presents the numerical results, and Section~\ref{section:conclusion} provides concluding remarks.

\section{Problem Formulation}\label{section:problem}
\subsection{Euler-Bernoulli Beam Mathematical Model}

We consider a flexible beam with a sliding boundary at $x=0$ and position and moment actuation at $x=1$. Without loss of generality, the beam length, mass density, and flexural rigidity are normalized to one. The equations are given as
\begin{align}
    & u_{tt}(t,x) + u_{xxxx}(t,x) = 0,\label{eq:euler_bernoulli_pde}\\
    & u_x(t,0) = u_{xxx}(t,0) = 0,\label{eq:euler_bernoulli_bc0}\\
    & u(t,1) = \theta_{1}(t), \quad u_{xx}(t,1) = \theta_{2}(t), \label{eq:euler_bernoulli_bc1}
\end{align}
where $x\in[0,1]$ is the space, $t\in[0,+\infty)$ is the time, $u$ is the displacement of the beam, and $\theta_{1}$ and $\theta_{2}$ are control inputs (position and moment actuation, respectively).

\subsection{Control Problem}
The goal of the ES method is to optimize an unknown static map $y = Q(\Theta)$ through real-time optimization, where $y^{*}$ and $\Theta^{*}$ denote the optimal unknown output and optimizer, respectively, while $y$ represents the measurable output, and $\theta_{1}$ and $\theta_{2}$ are the inputs.

In this work, the input of the map corresponds to the displacement at the uncontrolled end of the beam, which is subject to a sliding boundary condition (see Equation \eqref{eq:euler_bernoulli_bc0}). Thus, we define
\begin{align}
    \Theta(t) = u(t,0).
\end{align}

\begin{assum}[Local quadratic map]\label{assump:quadraticmap}
There exists $r_Q>0$ such that the map is exactly quadratic for $|\Theta-\Theta^*|<r_Q$:
\begin{equation}
    Q(\Theta(t)) = y^{*} + \frac{H}{2}(\Theta(t) - \Theta^{*})^{2},
    \label{extrachapter.eq:static_map}
\end{equation}
where $y^{*}, \Theta^{*} \in \mathbb{R}$, and $H < 0$ represents the Hessian.
\end{assum}

Assumption~\ref{assump:quadraticmap} is an exact local model; a general $C^2$ map includes a Taylor remainder, which is not covered by the present work.

\section{Extremum Seeking Boundary Control Design} \label{section:extremum_seeking_design}
\subsection{Demodulation Signals}
The demodulation signal $N(t)$, which is used to estimate the Hessian of the static map by multiplying it with the output $y(t)$ of the static map, is defined in \citet{Ghaffari2011MultivariableNE} as
\begin{align}
\hat{H}(t) = N(t)y(t),\quad &\text{with} \quad N(t) = -\dfrac{8}{a^{2}}\cos{(2\omega t)},
    \label{eq:hessian}
\end{align}
with $a\neq0$ and $\omega>0$, whereas the signal $M(t)$ is used to estimate the gradient of the static map as follows:
\begin{align}
    G(t) = M(t)y(t),\quad &\text{with} \quad M(t) = \frac{2}{a}\sin{(\omega t)}.\label{eq:gradient_estimate}
\end{align}

\subsection{Additive Probing Signal}
The boundary probing signals $S_1(t)$ and $S_2(t)$ are adapted from the basic ES to the case of PDE actuation dynamics \citep{Oliveira2017ExtremumSF,book2022}. The trajectory generation problem, as in \citet{Krsti2008BoundaryCO}, is described as follows:
\begin{align}
    R_{tt}(t,x) + R_{xxxx}(t,x)&=0, \label{eq:trajectoryPDE}\\
    R_{x}(t,0)=R_{xxx}(t,0)&=0,\label{eq:trajectorybc0}\\
     R_{xx}(t,0)&=0,\label{eq:trajectorybc02}\\
    R(t,0)&=a\sin(\omega t),\label{eq:trajectorybc03}\\
    S_{1}(t)&\coloneqq R(t,1),\label{eq:perturbation1}\\
    S_{2}(t)&\coloneqq R_{xx}(t,1).\label{eq:perturbation2}
\end{align}

The explicit solution of \eqref{eq:trajectoryPDE}-\eqref{eq:perturbation2} is stated in the following Lemma.

\begin{lem}[Boundary probing signals]
A smooth probing trajectory satisfying \eqref{eq:trajectoryPDE}--\eqref{eq:trajectorybc03} is
\begin{align}
R(t,x)=\frac a2[\cosh(\sqrt\omega x)+\cos(\sqrt\omega x)]\sin(\omega t).
\label{eq:solutionR}
\end{align}
The boundary probing signals in \eqref{eq:perturbation1}--\eqref{eq:perturbation2} are
\begin{align}
S_1(t)&=\frac a2[\cosh\sqrt\omega+\cos\sqrt\omega]\sin(\omega t),\label{eq:perturbationS1}\\
S_2(t)&=\frac{a\omega}{2}[\cosh\sqrt\omega-\cos\sqrt\omega]\sin(\omega t).\label{eq:perturbationS2}
\end{align}
\end{lem}
\begin{pf}
Write $R(t,x)$ for the probing trajectory of \eqref{eq:trajectoryPDE}-\eqref{eq:perturbation2} and consider its representation by the following power series:
    \begin{align}
        R(t,x)=\sum_{k=0}^{\infty}a_{k}(t)\frac{x^{k}}{k!}.\label{eq:series_trajectory_postulation}
    \end{align}
    
Substituting \eqref{eq:series_trajectory_postulation} into \eqref{eq:trajectoryPDE}-\eqref{eq:trajectorybc03}, it follows that
    \begin{equation}
    \begin{aligned}
    &a_{0}(t) = a\sin(\omega t) = a\mbox{Im}\{\mathrm{e}^{i\omega t}\},\\
        &a_{1}(t) = a_{2}(t)=a_{3}(t)=0, \\
        &a_{i+4}(t) = -\ddot{a}_{i}(t).
    \end{aligned}
    \label{eq:probe_coefficient_recursion}
    \end{equation}
    
Therefore,
    \begin{equation}
    \begin{aligned}
    &a_{4k}(t) = (-1)^{k}\frac{d^{2k}a_{0}}{dt^{2k}}(t) = \omega^{2k}a\sin(\omega t),\\
        &a_{4k+1}(t) = a_{4k+2}(t) = a_{4k+3}(t)=0.
    \end{aligned}
    \label{eq:probe_nonzero_coefficients}
    \end{equation}
    
The reference trajectory then becomes
    \begin{equation}
    \begin{aligned}
    R(t,x) &= \sum_{k=0}^{\infty}\omega^{2k}\frac{x^{4k}}{(4k)!}a\sin(\omega t), \\
        &= \frac{1}{2}\left[ \cosh (\sqrt{\omega}x)+\cos (\sqrt{\omega}x)\right]a\sin(\omega t).
    \end{aligned}
    \label{eq:probe_series_evaluation}
    \end{equation}
    
Finally, substituting \eqref{eq:solutionR} into \eqref{eq:perturbation1}-\eqref{eq:perturbation2}, respectively, we obtain \eqref{eq:perturbationS1}-\eqref{eq:perturbationS2}. The proof is complete.\hfill$\square$
\end{pf}

\subsection{Estimation Errors and PDE-Error Dynamics}
Since our objective is to drive the map input $\Theta(t)$ toward the unknown optimizer $\Theta^*$, we note that the corresponding nominal boundary values, after removing the probing signals, are $\hat{\theta}_1=\Theta^*$ and $\hat{\theta}_2=0$. We therefore define
\begin{align}
    \hat{\theta}_{1}(t) &\coloneqq \theta_{1}(t) - S_{1}(t),& \hat{\theta}_{2}(t) &\coloneqq \theta_{2}(t) - S_{2}(t), \label{eq:error_sigma} \\
    \hat{\Theta}(t) &\coloneqq \Theta(t) - a\sin(\omega t). \label{eq:error_theta1} 
\end{align}

Furthermore, we define the estimation errors in both the input and propagated input variables as
\begin{align}
    \tilde{\theta}_{1}(t) &\coloneqq \hat{\theta}_{1}(t) - \Theta^{*}, & \tilde{\theta}_{2}(t) &\coloneqq \hat{\theta}_{2}(t) - 0,
    \label{eq:error_thetatilde}\\
    \vartheta_{1}(t) &\coloneqq \hat{\Theta}(t) - \Theta^{*}, & \vartheta_{2}(t) &\coloneqq u_{xx}(t,0) - 0.\label{eq:error_vartheta}
\end{align}

Next, we define
\begin{equation}
    \alpha(t,x) \coloneqq u(t,x) - R(t,x) - \Theta^{*}. \label{eq:errorstate}
\end{equation}

Since $u$ and $R$ satisfy \eqref{eq:euler_bernoulli_pde} and \eqref{eq:trajectoryPDE}, respectively, and $\Theta^{*}$ is constant, the error defined in \eqref{eq:errorstate} satisfies
\begin{align}
    \alpha_{tt}+\alpha_{xxxx} = 0.\label{eq:errorPDE}
\end{align}

By differentiating \eqref{eq:errorstate} with respect to space once and three times, respectively, evaluating at $x=0$, and applying boundary conditions \eqref{eq:euler_bernoulli_bc0} and \eqref{eq:trajectorybc0}, we have
\begin{align}
    \alpha_{x}(t,0) = \alpha_{xxx}(t,0) = 0. \label{eq:errorPDE_bc}
\end{align}

Additionally, using $R(t,0)=a\sin(\omega t)$ and $R_{xx}(t,0)=0$, the definitions \eqref{eq:error_vartheta} and \eqref{eq:errorstate} yield
\begin{align}
 \vartheta_{1} = \alpha(t,0), & & \vartheta_{2}(t) = \alpha_{xx}(t,0). \label{eq:erroTheta}
\end{align}

Similarly, evaluating \eqref{eq:errorstate} and its second spatial derivative at $x=1$, and using \eqref{eq:euler_bernoulli_bc1}, \eqref{eq:perturbation1}--\eqref{eq:perturbation2}, and \eqref{eq:error_thetatilde}, we obtain
\begin{align}
    \alpha(t,1) = \tilde{\theta}_{1}(t), \qquad \alpha_{xx}(t,1) = \tilde{\theta}_{2}(t).\label{eq:errorbc1}
\end{align}

For sufficiently regular solutions with compatible initial and boundary data, define 
\begin{align}
  \beta(t,x)\coloneqq\alpha_t(t,x),  
\end{align}
and
\begin{equation}
\begin{aligned}
    U_1(t)&\coloneqq\dot{\tilde{\theta}}_1(t)
          =\dot{\hat{\theta}}_1(t),\\
    U_2(t)&\coloneqq\dot{\tilde{\theta}}_2(t)
          =\dot{\hat{\theta}}_2(t).
\end{aligned}
\label{eq:nominal_actuator_rates}
\end{equation}
Differentiating the error PDE \eqref{eq:errorPDE}, its boundary conditions \eqref{eq:errorPDE_bc} and \eqref{eq:errorbc1}, and \eqref{eq:erroTheta} with respect to time, we obtain the propagated error dynamics
\begin{align}
    &\dot{\vartheta}_{1}(t) = \beta(t,0),\label{eq:ODE1errordynamics}\\
    &\dot{\vartheta}_{2}(t) = \beta_{xx}(t,0),\label{eq:ODE2errordynamics}\\
    &\beta_{tt}(t,x)+\beta_{xxxx}(t,x) = 0,\label{eq:PDEerrordynamics}\\
    &\beta_{x}(t,0)=\beta_{xxx}(t,0) = 0,\label{eq:BC0errordynamics}\\
    &\beta(t,1) = U_{1}(t),\quad \beta_{xx}(t,1) = U_{2}(t),\label{eq:BCUerrordynamics}
\end{align}

\begin{figure*}
    \centering
    \includegraphics[width=0.75\linewidth]{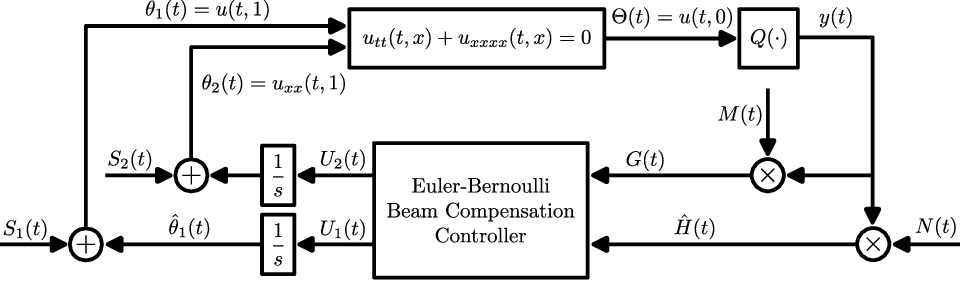}
    \caption{Block diagram of the boundary ES loop. The controller uses the PDE measurements in \eqref{eq:measured_feedback_fields} and includes the filter \eqref{eq:ES_filter_state}, the auxiliary state \eqref{eq:eta_extension}, and the integrator \eqref{eq:dynamic_control_U1} that generates $U_1$.}
    \label{fig:euler}
\end{figure*}

Adapting the proposed scheme in \citet{book2022} and combining \eqref{eq:euler_bernoulli_pde}-\eqref{extrachapter.eq:static_map}, the closed-loop ES with actuation dynamics governed by the EB beam PDE is illustrated in Figure \ref{fig:euler}.

\subsection{Euler--Bernoulli Beam Compensation via Backstepping Boundary Control}
We construct a backstepping boundary compensator using a Schr\"odinger-equation representation of the beam \citep{Smyshlyaev2009}, before incorporating the gradient and Hessian estimates into the feedback law.

For sufficiently regular solutions satisfying the initial and boundary compatibility conditions, introduce
\begin{align}
v(t,x)=\beta_t(t,x)-i\beta_{xx}(t,x).
\label{eq:transformation_schrodinger}
\end{align}
Using \eqref{eq:PDEerrordynamics}--\eqref{eq:BCUerrordynamics} gives the following Schr\"odinger-equation representation of the beam:
\begin{align}
v_t(t,x)+iv_{xx}(t,x)&=0,\label{eq:schrodingerPDE}\\
v(t,1)&=U(t),\label{eq:control_schrodinger}\\
U(t)&=\dot U_1(t)-iU_2(t),\label{eq:complex_control}\\
v_x(t,0)&=0.\label{eq:schrodinger_left_BC}
\end{align}

To express the coupling \eqref{eq:ODE1errordynamics}--\eqref{eq:ODE2errordynamics} at $x=0$ as a first-order system, introduce the auxiliary controller state
\begin{align}
\dot\eta(t)=\vartheta_2(t),\qquad\eta(0)=\eta_0,
\label{eq:eta_extension}
\end{align}
where $\eta_0\in\mathbb R$ is prescribed, and define
\begin{align}
X_1(t)&=\vartheta_1(t)-i\eta(t),\nonumber\\
X_2(t)&=\dot\vartheta_1(t)-i\vartheta_2(t)
=\beta(t,0)-i\vartheta_2(t).
\label{eq:complex_ODE_states}
\end{align}
We write $p(t)=\dot\vartheta_1(t)=\beta(t,0)$ for the real part of $X_2(t)$. Equations~\eqref{eq:ODE1errordynamics}, \eqref{eq:ODE2errordynamics}, and \eqref{eq:transformation_schrodinger} give
\begin{align}
\dot X(t)=AX(t)+Bv(t,0),
\label{eq:varthetadynamics}
\end{align}
where the column vectors $X=(X_1,\,X_2)$ and $B=(0,\,1)$ are written by components and $A=\left(\begin{smallmatrix}0&1\\0&0\end{smallmatrix}\right)$. The sliding-end conditions give the reconstruction
\begin{equation}
\begin{aligned}
\beta(t,x)&=\operatorname{Re}\{X_2(t)\}-\int_0^x(x-y)\operatorname{Im}\{v(t,y)\}\,dy,\\
\beta_t(t,x)&=\operatorname{Re}\{v(t,x)\}.
\end{aligned}
\label{eq:complex_beam_reconstruction}
\end{equation}

Equation~\eqref{eq:complex_beam_reconstruction} shows why $X_2$ is needed to recover the spatially constant component of $\beta$. Once the complex feedback $U$ is specified, the real control channels follow from \eqref{eq:complex_control}:
\begin{align}
\dot U_1(t)&=\operatorname{Re}\{U(t)\},\label{eq:real_control_channel1}\\
U_2(t)&=-\operatorname{Im}\{U(t)\}.\label{eq:real_control_channel2}
\end{align}

The channel \eqref{eq:real_control_channel1} is implemented as
\begin{align}
U_1(t)=U_1(0)+\int_0^t\operatorname{Re}\{U(\tau)\}\,d\tau.
\label{eq:dynamic_control_U1}
\end{align}

By \eqref{eq:nominal_actuator_rates}, the actuator variables are recovered from
\begin{align}
\hat\theta_1(t)&=\hat\theta_1(0)+\int_0^t U_1(\tau)\,d\tau,\label{eq:recovery_hattheta1}\\
\hat\theta_2(t)&=\hat\theta_2(0)+\int_0^t U_2(\tau)\,d\tau,\label{eq:recovery_hattheta2}
\end{align}
and the physical boundary inputs are
\begin{align}
\theta_1(t)&=\hat\theta_1(t)+S_1(t),\label{eq:physical_control_theta1}\\
\theta_2(t)&=\hat\theta_2(t)+S_2(t).\label{eq:physical_control_theta2}
\end{align}
Thus $\theta_1$ requires two integrations from $\operatorname{Re}\{U\}$, whereas $\theta_2$ requires one from $-\operatorname{Im}\{U\}$. Compatible classical data satisfy $U_1(0)=\beta(0,1)$ and $U_2(0)=\beta_{xx}(0,1)$.


\subsubsection{Target System}

We seek to map the system \eqref{eq:ODE1errordynamics}--\eqref{eq:BCUerrordynamics} into an exponentially stable cascade consisting of a two-dimensional complex ODE and a damped Schr\"odinger equation. The target system is
\begin{align}
\dot X(t) &= A_cX(t)+Bw(t,0),\label{eq:target_ode}\\
w_t(t,x)&=-iw_{xx}(t,x)-cw(t,x),\label{eq:target_pde}\\
w_x(t,0)&=w(t,1)=0,\label{eq:target_bc}
\end{align}
with
\begin{align}
A_c
&=A+BK_0
=
\begin{pmatrix}
0&1\\
-k_p&-k_d
\end{pmatrix},\\
K_0
&=
\begin{pmatrix}
-k_p&-k_d
\end{pmatrix}.
\label{eq:target_ode_matrices}
\end{align}

Here $K_0$ is the ODE feedback gain, $k_p,k_d>0$ specify the target dynamics, and $c>0$ is the beam damping parameter. The eigenvalues of \(A_c\) are
\begin{align}
\lambda_{1,2}=\frac{-k_d\pm\sqrt{k_d^2-4k_p}}{2}.
\label{eq:target_ode_eigenvalues}
\end{align}

For \(k_p>0\) and \(k_d>0\), both eigenvalues have strictly negative real parts. Therefore, \(A_c\) is Hurwitz.

As highlighted in \citet{Ren2013}, applying a single-step backstepping transformation to map \eqref{eq:schrodingerPDE}--\eqref{eq:schrodinger_left_BC} and \eqref{eq:varthetadynamics} into \eqref{eq:target_ode}--\eqref{eq:target_bc} is challenging due to the complexity of the associated kernels. To overcome this difficulty, a two-step design approach is adopted.

\textbf{First backstepping transformation:} Consider the ODE--Schr\"odinger system \eqref{eq:schrodingerPDE}--\eqref{eq:schrodinger_left_BC}, and \eqref{eq:varthetadynamics}. We introduce the backstepping transformation
\begin{align}
z(t,x)
&=
v(t,x)
-\int_0^xq(x,y)v(t,y)\,dy
-\gamma(x)X(t),
\label{eq:backsteppingtransformation1}
\end{align}
where \(\gamma(x)\in\mathbb C^{1\times2}\) and \(q(x,y)\in\mathbb C\).

The transformation \eqref{eq:backsteppingtransformation1} maps the ODE--Schr\"odinger system \eqref{eq:schrodingerPDE}--\eqref{eq:schrodinger_left_BC}, and \eqref{eq:varthetadynamics} into the intermediate target system
\begin{align}
\dot X(t)
&=
A_cX(t)+Bz(t,0),
\label{eq:intermediateODE}\\
z_t(t,x)
&=
-iz_{xx}(t,x),
\label{eq:intermediatePDE}\\
z_x(t,0)&=0,
\qquad
z(t,1)=W(t),
\label{eq:intermediateBC1}
\end{align}
where \(W(t)\in\mathbb C\) is an auxiliary boundary input to be determined in the second backstepping step.

Since \(\gamma(0)=K_0\), evaluating \eqref{eq:backsteppingtransformation1} at \(x=0\) gives \(z(t,0)=v(t,0)-K_0X(t)\), yielding \eqref{eq:intermediateODE}.

Differentiating \eqref{eq:backsteppingtransformation1} once with respect to time and twice with respect to space, and matching the resulting expression with \eqref{eq:intermediatePDE}--\eqref{eq:intermediateBC1}, gives the kernel equations
\begin{align}
\gamma''(x)
&=
i\gamma(x)A,
\label{eq:odekernel_intermediate}\\
\gamma'(0)&=0,\\
\gamma(0)&=K_0,
\label{eq:odekernelbc2_intermediate}\\
q_{yy}(x,y)-q_{xx}(x,y)
&=
0,
\label{eq:pdekernel_intermediate}\\
q_y(x,0)
&=
-i\gamma(x)B,\\
q(x,x)&=0.
\label{eq:pdekernelbc2_intermediate}
\end{align}

The solution of \eqref{eq:odekernel_intermediate}--\eqref{eq:odekernelbc2_intermediate} is
\begin{align}
\gamma(x)
=
\begin{pmatrix}
-k_p&
-k_d-i\dfrac{k_px^2}{2}
\end{pmatrix},
\qquad x\in[0,1].
\label{eq:solutionkernelode_intermediate}
\end{align}

The solution of \eqref{eq:pdekernel_intermediate}--\eqref{eq:pdekernelbc2_intermediate} is
\begin{align}
q(x,y)
&=
i\int_0^{x-y}\gamma(\sigma)B\,d\sigma
\nonumber\\
&=
\frac{k_p}{6}(x-y)^3
-ik_d(x-y).
\label{eq:solutionkernelpde_intermediate}
\end{align}

Finally, evaluating \eqref{eq:backsteppingtransformation1} at \(x=1\), using \(v(t,1)=U(t)\), and imposing \(z(t,1)=W(t)\), we obtain
\begin{align}
U(t)
&=
W(t)
+\int_0^1
\left[
\frac{k_p}{6}(1-y)^3
-ik_d(1-y)
\right]
v(t,y)\,dy
\nonumber\\
&\quad
-k_pX_1(t)
-\left(
k_d+i\frac{k_p}{2}
\right)X_2(t).
\label{eq:controller_schrodinger}
\end{align}
This is the complex-valued boundary control produced by the first backstepping transformation.

\textbf{Second backstepping transformation:} To obtain \eqref{eq:target_ode}--\eqref{eq:target_bc}, introduce
\begin{align}
w(t,x)
=
z(t,x)-\int_0^x\kappa(x,y)z(t,y)\,dy.
\label{eq:backsteppingtransformation2}
\end{align}
Since \(w(t,0)=z(t,0)\), the ODE in \eqref{eq:intermediateODE} becomes \eqref{eq:target_ode}. Differentiating \eqref{eq:backsteppingtransformation2}, substituting \eqref{eq:intermediatePDE}--\eqref{eq:intermediateBC1}, and matching \eqref{eq:target_pde}--\eqref{eq:target_bc} gives the kernel equations
\begin{align}
\kappa_{xx}(x,y)-\kappa_{yy}(x,y)
&=
ic\kappa(x,y),
\label{eq:kernelpdefinal}\\
\kappa_y(x,0)&=0,\\
\kappa(x,x)&=-i\frac{c}{2}x.
\label{eq:kernelbc2final}
\end{align}

The solution to the PDE \eqref{eq:kernelpdefinal}--\eqref{eq:kernelbc2final} is given on page 66 of \citet{Krsti2008BoundaryCO} by
\begin{align}
\kappa(x,y)
=
\kappa_r(x,y)+i\kappa_i(x,y),
\end{align}
\noindent where
\begin{align*}
\kappa_r(x,y)
&=
x\sqrt{\frac{c}{2(x^2-y^2)}}
\left[
-\mbox{ber}_1
\left(\sqrt{c(x^2-y^2)}\right)
\right.\\
&\left.
\hspace{2.5cm}
-\mbox{bei}_1
\left(\sqrt{c(x^2-y^2)}\right)
\right],\\
\kappa_i(x,y)
&=
x\sqrt{\frac{c}{2(x^2-y^2)}}
\left[
\mbox{ber}_1
\left(\sqrt{c(x^2-y^2)}\right)
\right.\\
&\left.
\hspace{2.5cm}
-\mbox{bei}_1
\left(\sqrt{c(x^2-y^2)}\right)
\right],
\end{align*}
and \(\mbox{ber}_1\) and \(\mbox{bei}_1\) are the Kelvin functions of order one.

Then, evaluating \eqref{eq:backsteppingtransformation2} at \(x=1\) and using \eqref{eq:target_bc}, we obtain
\begin{align}
W(t)
=
\int_0^1\kappa(1,y)z(t,y)\,dy.
\label{eq:intermediatecontroller}
\end{align}

Substituting \eqref{eq:backsteppingtransformation1} into \eqref{eq:intermediatecontroller}, changing the order of integration, and using \eqref{eq:controller_schrodinger} gives the composed feedback law
\begin{align}
U(t)
=
\int_0^1\ell(y)v(t,y)\,dy+hX(t),
\label{eq:complex_control_law_mean}
\end{align}
where
\begin{equation}
\begin{aligned}
\ell(y)
&=
q(1,y)+\kappa(1,y)
-\int_y^1\kappa(1,\xi)q(\xi,y)\,d\xi,\\
h
&=
\gamma(1)
-\int_0^1\kappa(1,y)\gamma(y)\,dy.
\end{aligned}
\label{eq:composed_boundary_gains}
\end{equation}


\subsubsection{Invertibility of the Transformations}
Following the Volterra-operator arguments in \citet{Ren2013}, the transformations \eqref{eq:backsteppingtransformation1} and \eqref{eq:backsteppingtransformation2} are invertible. Specifically, the inverse of \eqref{eq:backsteppingtransformation1} can be written as
\begin{equation}
\begin{aligned}
v(t,x)
=
z(t,x)
+\int_0^x\iota(x,y)z(t,y)\,dy
+\psi(x)X(t),
\end{aligned}
\label{eq:inverse_first_step}
\end{equation}
where \(\iota(x,y)\in\mathbb C\) and \(\psi(x)\in\mathbb C^{1\times2}\). Similarly, the inverse of \eqref{eq:backsteppingtransformation2} is given by
\begin{equation}
\begin{aligned}
z(t,x)
=
w(t,x)
+\int_0^x\rho(x,y)w(t,y)\,dy.
\end{aligned}
\label{eq:inverse_second_step}
\end{equation}
The inverse kernels are uniquely specified by
\begin{align}
\iota(x,y)&=q(x,y)+\int_y^xq(x,s)\iota(s,y)\,ds,\nonumber\\
\rho(x,y)&=\kappa(x,y)+\int_y^x\kappa(x,s)\rho(s,y)\,ds,\nonumber\\
\psi(x)&=\gamma(x)+\int_0^x\iota(x,y)\gamma(y)\,dy.
\label{eq:inverse_Volterra_resolvents}
\end{align}
{\color{blue}By} \citet[Lemmas~1--2]{Ren2013}, the transformations and their inverses are bounded on the corresponding ODE--PDE state spaces. The smoothness of the kernels also ensures that these transformations preserve $H^2(0,1)$ regularity.

\subsection{Target System of the Euler--Bernoulli Beam Equation}

To determine the target system in the original Euler--Bernoulli coordinates, we use transformations \eqref{eq:backsteppingtransformation1} and \eqref{eq:backsteppingtransformation2}, together with the control law \eqref{eq:complex_control_law_mean}, and define
\begin{align}
\zeta(t,x) = \int_x^1\int_0^y\operatorname{Im}\{w(t,\xi)\}\,d\xi\,dy.
\label{eq:transformation_bernoulli}
\end{align}
Here, \(w\) is the state of the Schr\"odinger target system \eqref{eq:target_pde}--\eqref{eq:target_bc}. Computing the second-order time derivative and the fourth-order spatial derivative of \eqref{eq:transformation_bernoulli}, integrating by parts, and using the boundary conditions in \eqref{eq:target_bc}, we verify that \(\zeta\) satisfies \citep{Smyshlyaev2009}
\begin{align}
\zeta_{tt}(t,x)+2c\zeta_t(t,x) +c^2\zeta(t,x)+\zeta_{xxxx}(t,x)&=0,
\label{eq:target_bernoulli_pde}\\
\zeta_x(t,0)= \zeta_{xxx}(t,0) &=0, \label{eq:target_bernoulli_bc_left}\\
\zeta(t,1) = \zeta_{xx}(t,1) &= 0. \label{eq:target_bernoulli_bc4}
\end{align}

The finite-dimensional component follows from \eqref{eq:target_ode}, since both backstepping transformations leave $X$ unchanged. From \eqref{eq:transformation_bernoulli} and \eqref{eq:target_pde}--\eqref{eq:target_bc}, we obtain
\begin{equation}
w(t,x)=\zeta_t(t,x)+c\zeta(t,x)-i\zeta_{xx}(t,x).
\label{eq:target_beam_complex_relation}
\end{equation}
Although \eqref{eq:ODE1errordynamics} is written as $\dot\vartheta_1(t)=\beta(t,0)$, the boundary error velocity $\beta(t,0)=\operatorname{Re}\{X_2(t)\}$ also evolves dynamically. Consequently, the real and negative imaginary parts of the equation for $\dot X_2$ in \eqref{eq:target_ode} give $\ddot\vartheta_1$ and $\dot\vartheta_2$, respectively. Using \eqref{eq:target_beam_complex_relation}, we obtain
\begin{align}
\ddot{\vartheta}_1(t)
&=-k_p\vartheta_1(t)-k_d\dot{\vartheta}_1(t)\nonumber\\
&\quad+\zeta_t(t,0)+c\zeta(t,0),
\label{eq:target_beam_vartheta1}\\
\dot{\vartheta}_2(t)
&=-k_p\eta(t)-k_d\vartheta_2(t)+\zeta_{xx}(t,0),
\label{eq:target_beam_vartheta2}
\end{align}
where $\dot\eta=\vartheta_2$ by \eqref{eq:eta_extension}. Thus, the complete target system consists of the damped beam \eqref{eq:target_bernoulli_pde}--\eqref{eq:target_bernoulli_bc4} driving the finite-dimensional subsystem \eqref{eq:target_beam_vartheta1}--\eqref{eq:target_beam_vartheta2}.

\subsection{Control Laws}
The feedback \eqref{eq:complex_control_law_mean} can be expressed directly in the Euler--Bernoulli error coordinates. Write $\ell=\ell_r+i\ell_i$ and $h=(h_1,h_2)$, with $h_j=h_{j,r}+ih_{j,i}$, $j\in\{1,\,2\}$. Using \eqref{eq:transformation_schrodinger}, \eqref{eq:complex_ODE_states}, and \eqref{eq:complex_control}, we obtain
\begin{align}
\dot U_1(t)
&=\int_0^1\Bigl[\ell_r(y)\beta_t(t,y)
+\ell_i(y)\beta_{xx}(t,y)\Bigr]\,dy
\nonumber\\
&\quad+h_{1,r}\vartheta_1(t)+h_{1,i}\eta(t)
\nonumber\\
&\quad+h_{2,r}\beta(t,0)+h_{2,i}\vartheta_2(t),
\label{eq:control_law_beam_deterministic1}\\
U_2(t)
&=\int_0^1\Bigl[\ell_r(y)\beta_{xx}(t,y)
-\ell_i(y)\beta_t(t,y)\Bigr]\,dy
\nonumber\\
&\quad-h_{1,i}\vartheta_1(t)+h_{1,r}\eta(t)
\nonumber\\
&\quad-h_{2,i}\beta(t,0)+h_{2,r}\vartheta_2(t).
\label{eq:control_law_beam_deterministic2}
\end{align}

The auxiliary state $\eta$ satisfies \eqref{eq:eta_extension}, and $U_1$ is obtained by integration as in \eqref{eq:dynamic_control_U1}. The physical boundary inputs follow from \eqref{eq:recovery_hattheta1}--\eqref{eq:physical_control_theta2}.

\subsubsection{Implementable Extremum Seeking Control Law}

The feedback laws \eqref{eq:control_law_beam_deterministic1}--\eqref{eq:control_law_beam_deterministic2} cannot be implemented directly because $\vartheta_1$ depends on the unknown optimizer $\Theta^*$. In this section, we construct an implementable realization using the gradient and Hessian estimates.

Choose fixed gains $\mu_p,\mu_d>0$ independently of $H$. For fixed $\vartheta_1$ with $|\vartheta_1|+|a|<r_Q$, let $\langle\cdot\rangle$ denote the average over one probing period with the state held fixed. Then \citep[Section~3.2.1]{OliveiraKrstic2023}
\begin{align}
\langle G\rangle=H\vartheta_1,\qquad \langle\hat H\rangle=H.
\label{eq:hessian_estimates}
\end{align}
To relate these gains to the design in \eqref{eq:target_ode_matrices}, set $k_p=-H\mu_p>0$ and $k_d=-H\mu_d>0$ for analysis only. Define the coefficients independent of $H$:
\begin{equation}
\begin{aligned}
\bar q(x,y)&=-\frac{\mu_p}{6}(x-y)^3+i\mu_d(x-y),\\
\bar\gamma(x)&=\begin{pmatrix}\mu_p&\mu_d+i\mu_px^2/2\end{pmatrix},\\
\bar\ell(y)&=\bar q(1,y)-\int_y^1\kappa(1,\xi)\bar q(\xi,y)\,d\xi,\\
\bar h&=\bar\gamma(1)-\int_0^1\kappa(1,y)\bar\gamma(y)\,dy.
\end{aligned}
\label{eq:normalized_feedback_gains}
\end{equation}
Then \eqref{eq:solutionkernelode_intermediate}--\eqref{eq:solutionkernelpde_intermediate} give $\gamma=H\bar\gamma$ and $q=H\bar q$, so \eqref{eq:composed_boundary_gains} yields
\begin{align}
\ell(y)=\kappa(1,y)+H\bar\ell(y),\qquad h=H\bar h.
\label{eq:feedback_split}
\end{align}

Using \eqref{eq:complex_ODE_states}, we replace $H\vartheta_1$ by $G$ and the remaining factors $H$ by $\hat H$. The known term involving $\kappa$ is applied directly to retain the damping in the PDE system, while the map-dependent terms are low-pass filtered. Writing $\bar h=(\bar h_1,\, \bar h_2)$ and $p(t)=\beta(t,0)$, the filter state $F(t)\in\mathbb C$ and the boundary control satisfy
\begin{align}
\dot F
&=-\overline c F+\overline c\bar h_1G
+\overline c\hat H\Bigl[\int_0^1\bar\ell(y)v(t,y)\,dy
\nonumber\\
&\qquad-i\bar h_1\eta+\bar h_2(p-i\vartheta_2)\Bigr],
\label{eq:ES_filter_state}\\
U(t)
&=\int_0^1\kappa(1,y)v(t,y)\,dy+F(t),
\qquad \overline c>0.
\label{eq:split_filtered_complex_control}
\end{align}
The corresponding control laws in beam coordinates are
\begin{align}
\dot U_1(t)
&=\int_0^1\kappa_r(1,y)\beta_t(t,y)\,dy
\nonumber\\
&\quad+\int_0^1\kappa_i(1,y)\beta_{xx}(t,y)\,dy
+\operatorname{Re}\{F(t)\},
\label{eq:control_law_filtered}\\
U_2(t)
&=-\int_0^1\kappa_i(1,y)\beta_t(t,y)\,dy
\nonumber\\
&\quad+\int_0^1\kappa_r(1,y)\beta_{xx}(t,y)\,dy
-\operatorname{Im}\{F(t)\}.
\label{eq:control_law_filtered_2}
\end{align}

Implementation of \eqref{eq:control_law_filtered}-\eqref{eq:control_law_filtered_2} requires the measured output $y$, the boundary traces $u_t(t,0)$ and $u_{xx}(t,0)$, and the distributed derivatives $u_{tt}$ and $u_{txx}$. By \eqref{eq:solutionR} and \eqref{eq:errorstate},
\begin{equation}
\begin{aligned}
p(t)&=u_t(t,0)-a\omega\cos(\omega t),\\
\vartheta_2(t)&=u_{xx}(t,0),\\
\beta_t&=u_{tt}-R_{tt},\qquad
\beta_{xx}=u_{txx}-R_{txx}.
\end{aligned}
\label{eq:measured_feedback_fields}
\end{equation}
Thus, evaluating the controller requires no knowledge of $\Theta^*$, $y^*$, or $H$. The auxiliary state and actuator realization follow \eqref{eq:eta_extension} and \eqref{eq:dynamic_control_U1}--\eqref{eq:physical_control_theta2}.

For classical solutions, initialize $\hat\theta_1(0)=u(0,1)-S_1(0)$ and $\hat\theta_2(0)=u_{xx}(0,1)-S_2(0)$, impose the compatibility conditions stated after \eqref{eq:physical_control_theta2}, and select
\begin{align}
F(0)=v(0,1)-\int_0^1\kappa(1,y)v(0,y)\,dy.
\label{eq:ES_filter_initial_compatibility}
\end{align}

\section{Stability Analysis of the Closed-Loop System}\label{section:stability}

We first prove exponential stability of the origin of the averaged closed-loop system with the filter dynamics included. For fixed $0<|a|<r_Q$, we then use a high-frequency estimate to prove that the full ES error system admits a locally exponentially attracting periodic solution for sufficiently large $\omega$.

Consider the error dynamics \eqref{eq:errorPDE}--\eqref{eq:BCUerrordynamics} under the boundary feedback \eqref{eq:control_law_filtered}--\eqref{eq:control_law_filtered_2} and the filter \eqref{eq:ES_filter_state}. With $f=(\operatorname{Re}\{F\},\,-\operatorname{Im}\{F\})$, the closed-loop error state $Z=(\alpha,\,\beta,\,\eta,\, f)$ belongs to
\begin{equation}
\begin{aligned}
\mathcal X_{\rm EB}
=\bigl\{Z &\in H^4\times H^2\times\mathbb R\times\mathbb R^2 : {}\\
&\alpha_x(0)=\alpha_{xxx}(0)=\beta_x(0)=0\bigr\}.
\end{aligned}
\label{eq:physical_EB_state_space}
\end{equation}
We equip this space with the norm
\begin{align}
\mathcal N_{\rm EB}(Z)^2
=\|\alpha\|_{H^4}^2+\|\beta\|_{H^2}^2+\eta^2+|f|^2,
\label{eq:physical_EB_state_norm}
\end{align}
where $H^m=H^m(0,1)$ and $\|h\|_{H^m}^2=\sum_{j=0}^m\|h^{(j)}\|_{L^2}^2$. Since $\beta_t=-\alpha_{xxxx}$, both $\beta_t$ and $\beta_{xx}$ belong to $L^2(0,1)$, as required by the boundary feedback. We derive the estimates for classical solutions with compatible initial data. Appendix~\ref{app:beam_coordinates} gives the generator domain and extends the estimates to mild solutions.

\subsection{Averaged Closed-Loop System and Backstepping Transformation}\label{sec:beam_filter_dynamics}

We average over the probing phase with the Euler--Bernoulli error state $Z$ held fixed. This also fixes $v$ and $X$ through \eqref{eq:transformation_schrodinger} and \eqref{eq:complex_ODE_states}, since $\beta_t=-\alpha_{xxxx}$. We use the Schr\"odinger representation to combine the two real feedback channels in one equation and identify the boundary input to the backstepping target system. By \eqref{eq:feedback_split}, the map-dependent term in the control law \eqref{eq:complex_control_law_mean} is $H\Phi$, where
\begin{align}
\Phi(t)=\int_0^1\bar\ell(y)v(t,y)\,dy+\bar hX(t).
\label{eq:map_dependent_feedback_signal}
\end{align}
Since $X_1=\vartheta_1-i\eta$, the expression multiplied by $\hat H$ in \eqref{eq:ES_filter_state} is $\Phi-\bar h_1\vartheta_1$. Adding and subtracting $H\Phi$ in \eqref{eq:ES_filter_state} gives
\begin{align}
\dot F
&=-\overline c(F-H\Phi)+\overline c(\hat H-H)\Phi\nonumber\\
&\quad+\overline c\bar h_1(G-\hat H\vartheta_1).
\label{eq:exact_complex_closed_loop_filter}
\end{align}

For $|\vartheta_1|+|a|<r_Q$, \eqref{eq:hessian_estimates} gives the phase averages $H\vartheta_1$ and $H$ of $G$ and $\hat H$, respectively. The last two terms in \eqref{eq:exact_complex_closed_loop_filter} therefore have zero phase average. Since the probing phase enters the error dynamics only through these terms, the averaged closed loop is obtained by replacing the filter equation with
\begin{align}
\dot F_{\rm av}=-\overline c(F_{\rm av}-H\Phi_{\rm av}).
\label{eq:averaged_complex_filter}
\end{align}

Proposition~\ref{prop:filtered_average_stability} analyzes the linear extension of these averaged equations to $\mathcal X_{\rm EB}$. Within $|\vartheta_1|+|a|<r_Q$, they coincide with the phase average of the ES error system.

We now omit the subscript $\mathrm{av}$ and apply the backstepping transformations \eqref{eq:backsteppingtransformation1} and \eqref{eq:backsteppingtransformation2} with $k_p=-H\mu_p$ and $k_d=-H\mu_d$. This gives the target ODE \eqref{eq:target_ode} and PDE \eqref{eq:target_pde}. At $x=1$, the transformations subtract the deterministic feedback \eqref{eq:complex_control_law_mean} from the applied input \eqref{eq:split_filtered_complex_control}. The integral involving $\kappa$ cancels, leaving the filter error as the boundary input to the target PDE. The boundary conditions and filter dynamics are therefore
\begin{align}
w_x(t,0)&=0,\qquad w(t,1)=e_f(t),\nonumber\\
e_f(t)&=F(t)-H\Phi(t),\qquad \dot F(t)=-\overline c e_f(t).
\label{eq:filtered_target_boundary}
\end{align}
The inverse transformations \eqref{eq:inverse_second_step} and \eqref{eq:inverse_first_step} recover $v$ from $(X,w)$ and hence determine $\Phi$ in \eqref{eq:map_dependent_feedback_signal}.

\subsection{Exponential Stability of the Averaged Closed-Loop System}

To prove decay in the norm \eqref{eq:physical_EB_state_norm}, we use the coordinates $(z_1,z_2,z_3,z_4,f,\chi,\upsilon)$ defined by \eqref{eq:real_beam_predictors} and \eqref{eq:real_beam_lifting}. Identity~\eqref{eq:full_target_energy_coordinate_bridge} relates these coordinates to $w$. The variables $(\chi,\upsilon)$ satisfy homogeneous boundary conditions, which we use in the energy estimate. The filter in~\eqref{eq:real_predictor_feedback_identity} and the four scalar beam coordinates in~\eqref{eq:real_beam_predictors} form a finite-dimensional subsystem. Suppressing the subscript $\mathrm{av}$, set

\begin{align}
\xi&=(z_1,z_2,f_1,z_3,z_4,f_2)^T,\nonumber\\
M_3&=\begin{pmatrix}0&1&0\\0&0&1\\
\overline cH\mu_p&\overline cH\mu_d&-\overline c\end{pmatrix},\nonumber\\  M_f&=\operatorname{diag}(M_3,M_3).
\label{eq:real_filter_matrix}
\end{align}

Equations~\eqref{eq:real_predictor_dynamics} and \eqref{eq:real_predictor_feedback_identity} give $\dot\xi=M_f\xi$. Since $M_f=\operatorname{diag}(M_3,M_3)$, it suffices to examine the characteristic polynomial of $M_3$,
\begin{align}
\det(sI-M_3)=s^3+\overline c s^2-\overline cH\mu_ds-\overline cH\mu_p.
\label{eq:real_filter_polynomial}
\end{align}
All coefficients are positive because $H<0$ and $\overline c,\mu_p,\mu_d>0$. The remaining Routh--Hurwitz inequality for this cubic is
$\overline c(-\overline cH\mu_d)>-\overline cH\mu_p$.
Dividing by $-\overline cH>0$ shows that the finite-dimensional subsystem is exponentially stable if and only if
\begin{align}
\overline c>\mu_p/\mu_d.
\label{eq:filter_cutoff_condition}
\end{align}
The filter bandwidth \eqref{eq:filter_cutoff_condition} can therefore be chosen using the design gains alone.

To establish exponential stability of the full averaged system, it remains to estimate the response of $(\chi,\upsilon)$ to $\dot f$ in~\eqref{eq:lifted_real_beam_dynamics}. The transformation from $Z=(\alpha,\beta,\eta,f)$ to $(\xi,\chi,\upsilon)$ uses \eqref{eq:real_beam_predictors} to define $z_1,z_2,z_3$, and $z_4$, \eqref{eq:real_filter_matrix} to collect these coordinates and $f$ in $\xi$, and \eqref{eq:real_beam_lifting} to define $(\chi,\upsilon)$. Appendix~\ref{app:beam_coordinates} proves that this transformation and its inverse are bounded and gives the constants $m_1,M_1>0$ in~\eqref{eq:explicit_coordinate_constants} for which
\begin{align}
m_1\mathcal N_{\rm EB}(Z)^2
&\leq |\xi|^2+\|\chi_{xx}\|_{L^2}^2+\|\upsilon\|_{L^2}^2\nonumber\\
&\leq M_1\mathcal N_{\rm EB}(Z)^2.
\label{eq:real_beam_final_norm_equivalence}
\end{align}
Thus an exponential decay estimate in the transformed coordinates yields the same decay rate in the original error norm~\eqref{eq:physical_EB_state_norm}.

The following proposition shows that the bandwidth condition
\eqref{eq:filter_cutoff_condition}, derived for the finite-dimensional
subsystem, also ensures exponential stability of the full averaged
closed-loop system in the norm~\eqref{eq:physical_EB_state_norm}.

\begingroup\wrapresultname
\begin{prop}[Averaged closed loop is exponentially stable for $\overline c>\mu_p/\mu_d$]\label{prop:filtered_average_stability}
Consider the averaged version of the EB error system \eqref{eq:errorPDE}--\eqref{eq:BCUerrordynamics} with the auxiliary state~\eqref{eq:eta_extension}, the boundary control laws \eqref{eq:control_law_filtered}--\eqref{eq:control_law_filtered_2}, and the averaged filter~\eqref{eq:averaged_complex_filter}. Let $Z_{\rm av}=(\alpha_{\rm av},\beta_{\rm av},\eta_{\rm av},f_{\rm av})$ denote the averaged closed-loop error state vector. If $H<0$, $c,\mu_p,\mu_d>0$, and $\overline c>\mu_p/\mu_d$, then the averaged closed-loop system is
well posed on $\mathcal X_{\rm EB}$. For every initial state $Z_{\rm av}(0)\in\mathcal X_{\rm EB}$, its mild solution satisfies
\begin{align}
\mathcal N_{\rm EB}(Z_{\rm av}(t))
&\leq C_{\rm EB}e^{-\lambda t}\mathcal N_{\rm EB}(Z_{\rm av}(0)),
\quad t\geq0,
\label{eq:physical_EB_exponential_stability}
\end{align}
with
\begin{align}
C_{\rm EB}=\sqrt{\frac{M_1q_+}{m_1q_-}},\qquad
\lambda=\frac{\sigma_f}{2}{\color{blue},}
\label{eq:explicit_average_constants}
\end{align}
where $m_1,M_1$ are given in~\eqref{eq:real_beam_final_norm_equivalence} and $Q=Q^T>0$ solves $M_f^TQ+QM_f=-I_6$ with $M_f$ from~\eqref{eq:real_filter_matrix}. For $k=\sqrt{c/2}$, define
\begin{equation}
D_c^2=\frac{\sinh(2k)+\sin(2k)}{2k[\cosh(2k)+\cos(2k)]}
\label{eq:average_lifting_constant}
\end{equation}
and set
\begin{equation}
\begin{aligned}
\rho &=1+\frac{D_c^2\overline c^2}{2c}
             [1+H^2(\mu_p^2+\mu_d^2)],\\
\sigma_f&=\min\{[\rho\lambda_{\max}(Q)]^{-1},c\},\\
q_-&=\min\{\rho\lambda_{\min}(Q),1/2\},\\
q_+&=\max\{\rho\lambda_{\max}(Q),1/2\}.
\end{aligned}
\label{eq:average_decay_coefficients}
\end{equation}
\end{prop}
\endgroup
\begin{pf}
Define $J_f\in\mathbb R^{2\times6}$ by $J_f\xi=(f_1,f_2)^T$ and set $L_f=J_fM_f$. Since $\dot\xi=M_f\xi$, the input to the PDE \eqref{eq:lifted_real_beam_dynamics} is $\dot f=L_f\xi$.
Under \eqref{eq:filter_cutoff_condition}, $M_3$ is Hurwitz, and hence so is $M_f=\operatorname{diag}(M_3,M_3)$. Let $Q=Q^T>0$ be the unique solution in $\mathbb R^{6\times6}$ of $M_f^TQ+QM_f=-I_6$. Consequently, along every solution of the averaged closed-loop system, $\xi$ converges exponentially to zero.

Consider the Lyapunov functional
\begin{align}
V_f=\rho\xi^TQ\xi+E_b,\label{eq:lyap_average}
\end{align}
where
\begin{equation}
E_b=\tfrac12\bigl(\|\chi_{xx}\|_{L^2}^2+\|\upsilon\|_{L^2}^2\bigr).
\label{eq:average_lifted_energy}
\end{equation}
For compatible classical solutions, differentiate \eqref{eq:lyap_average} using $\dot\xi=M_f\xi$ from~\eqref{eq:real_lifted_average_system} and the PDE~\eqref{eq:lifted_real_beam_dynamics}. Integrating by parts in $x$ and applying the boundary conditions~\eqref{eq:lifted_real_beam_boundaries} gives
\begin{equation}
\begin{aligned}
\dot V_f
&=-\rho|\xi|^2-2cE_b\\
&\quad+\int_0^1\chi_{xx}(b_c\dot f_1-a_c\dot f_2)\,dx\\
&\quad-\int_0^1\upsilon(a_c\dot f_1+b_c\dot f_2)\,dx.
\end{aligned}
\label{eq:average_energy_derivative}
\end{equation}
Since $\dot f=L_f\xi$ and $\|L_f\|^2=\overline c^2[1+H^2(\mu_p^2+\mu_d^2)]$, the Cauchy--Schwarz inequality gives
\begin{equation}
\dot V_f\leq-\rho|\xi|^2-2cE_b+D_c|\dot f|\sqrt{2E_b}{\color{blue},}
\label{eq:average_cauchy_schwarz}
\end{equation}
with $\rho$ from Proposition~\ref{prop:filtered_average_stability}. The formula for $D_c^2=\int_0^1(a_c^2+b_c^2)\,dx$ follows from~\eqref{eq:real_lifting_explicit}. Young's inequality gives
\begin{equation}
D_c|\dot f|\sqrt{2E_b}
\leq cE_b+\frac{D_c^2\|L_f\|^2}{2c}|\xi|^2.
\label{eq:average_young_bound}
\end{equation}
Thus, with $\rho$ defined in \eqref{eq:average_decay_coefficients}, it follows that
\begin{align}
\dot V_f\leq-|\xi|^2-cE_b\leq-\sigma_fV_f{\color{blue}.}
\label{eq:filtered_beam_Lyapunov_decay}
\end{align}
Set
\begin{equation}
m_V=m_1q_-,\qquad M_V=M_1q_+.
\label{eq:average_energy_norm_constants}
\end{equation}
Equation~\eqref{eq:real_beam_final_norm_equivalence} gives
\begin{align}
m_V\mathcal N_{\rm EB}(Z_{\rm av}(t))^2
\leq V_f(t)\leq M_V\mathcal N_{\rm EB}(Z_{\rm av}(t))^2.
\label{eq:Lyapunov_physical_norm_bounds}
\end{align}
Integration of \eqref{eq:filtered_beam_Lyapunov_decay} and application of \eqref{eq:Lyapunov_physical_norm_bounds} give \eqref{eq:physical_EB_exponential_stability} for compatible classical solutions with the constants in \eqref{eq:explicit_average_constants}.

Appendix~\ref{app:beam_coordinates} establishes well-posedness on $\mathcal X_{\rm EB}$ and extends this estimate to mild solutions. \hfill $\square$
\end{pf}


\subsection{Convergence of the Extremum-Seeking Closed-Loop}

Let $Y=(\xi,\chi,\upsilon)$, with $\|Y\|^2=|\xi|^2+\|\chi_{xx}\|_{L^2}^2+\|\upsilon\|_{L^2}^2$, and let $\mathcal L$ be the averaged generator defined in Appendix~\ref{app:beam_coordinates}. Proposition \ref{prop:filtered_average_stability}  gives $\|e^{\mathcal Lt}\|\leq M_Le^{-\lambda_Lt}$ with
\begin{align}
M_L=\sqrt{q_+/q_-},\qquad \lambda_L=\sigma_f/2.
\label{eq:explicit_semigroup_constants}
\end{align}

Collect the coordinates entering the ES correction as
\begin{align}
\mathbf r(Y)=(z_1,z_2,z_3,z_4,\vartheta_1)^T.
\label{eq:ES_error_vector}
\end{align}
Using \eqref{eq:exact_predictor_ES_signal}, the exact dynamics can be written as
\begin{align}
\dot Y=\mathcal LY+\mathcal B\Delta(\omega t,\mathbf r(Y)).
\label{eq:real_semilinear_ES_system}
\end{align}
Appendix~\ref{app:beam_averaging} gives the bounded operator $\mathcal B$ and the explicit smooth perturbation $\Delta$, which has zero phase mean and only four harmonics. The boundedness of $\mathbf r$ follows from the reconstruction in Appendix~\ref{app:beam_coordinates}.

\begin{rem}The averaging result for analytic semigroups in \citet{OliveiraKrstic2023} cannot be applied directly: the beam modes have fixed decay rate $c$ and unbounded oscillation frequencies, giving a nonanalytic $C_0$-semigroup. We instead estimate its frequency response along the forcing directions and construct the periodic correction explicitly in Appendix~\ref{app:beam_averaging}.
\end{rem}

\begingroup\wrapresultname
\begin{thm}[ES error converges to an $\mathcal O(\omega^{-1/2})$ periodic orbit]\label{thm:ES_stability_revised}
Consider the Euler--Bernoulli PDE \eqref{eq:euler_bernoulli_pde}--\eqref{eq:euler_bernoulli_bc1}
with boundary inputs \eqref{eq:physical_control_theta1}--\eqref{eq:physical_control_theta2}
generated by the control laws \eqref{eq:control_law_filtered}--\eqref{eq:control_law_filtered_2}
and the filter \eqref{eq:ES_filter_state}, equivalently \eqref{eq:exact_predictor_ES_signal} in real coordinates. Suppose Assumption~\ref{assump:quadraticmap} holds and fix $0<|a|<r_Q$, $c,\mu_p,\mu_d>0$, and $\overline c>\mu_p/\mu_d$. There exist $r_*,\omega_0,C,\delta,\lambda_*,M_*>0$ such that, for every $\omega\geq\omega_0$ and $\Pi=2\pi/\omega$, the closed-loop error system \eqref{eq:errorPDE}--\eqref{eq:BCUerrordynamics} admits a $\Pi$-periodic solution $Z^\Pi$. It is unique among $\Pi$-periodic error solutions satisfying $\sup_{t\in\mathbb R}\mathcal N_{\rm EB}(Z(t))<r_*/(2\sqrt{M_1})$.

For every $Z(0)\in\mathcal X_{\rm EB}$ satisfying $\mathcal N_{\rm EB}(Z(0)-Z^\Pi(0))<\delta$, the mild solution exists for all $t\geq0$, satisfies $|\Theta(t)-\Theta^*|<r_Q$, and obeys
\begin{align}
\mathcal N_{\rm EB}(Z(t)-Z^\Pi(t))
&\leq M_*e^{-\lambda_*t}
\mathcal N_{\rm EB}(Z(0)-Z^\Pi(0)).
\label{eq:uniform_physical_orbit_attraction}
\end{align}
The periodic error satisfies
\begin{align}
\sup_{t\in\mathbb R}\mathcal N_{\rm EB}(Z^\Pi(t))
\leq C\omega^{-1/2}.
\label{eq:convergance_closedloop}
\end{align}
For the same initial data,
\begin{align}
\limsup_{t\to\infty}|\Theta(t)-\Theta^*|
&\leq |a|+C\omega^{-1/2},\label{eq:convergence_Theta}\\
\limsup_{t\to\infty}|\theta_1(t)-\Theta^*|
&\leq r_1(\omega)+C\omega^{-1/2},\label{eq:convergence_theta1}\\
\limsup_{t\to\infty}|y(t)-y^*|
&\leq\tfrac{|H|}{2}(|a|+C\omega^{-1/2})^2,
\label{eq:convergence_y}
\end{align}
where $r_1(\omega)=|a|\,|\cosh\sqrt\omega+\cos\sqrt\omega|/2$ is the amplitude of the position probe in \eqref{eq:perturbationS1}. The constants $r_*,\omega_0,C,\delta,\lambda_*,M_*$ depend only on $H,y^*,r_Q,a,\allowbreak c,\overline c,\mu_p,\mu_d$ and are uniform in $\omega\geq\omega_0$. With $d=\sqrt{1+2/m_1}$, $r_*=(r_Q-|a|)/(4d)$, and $C_0=M_LK/\lambda_L$, one choice is
\begin{align}
\omega_0&=\left[\max\left\{1,2A_1,
\frac{4\sqrt{M_1/m_1}(A_0+C_0)}{r_*},2C_0\right\}\right]^2,\nonumber\\
\delta&=\frac{r_*}{8M_L\sqrt{M_1}},\qquad
\lambda_*=\frac{\lambda_L}{2},\nonumber\\
M_*&=3M_L\sqrt{M_1/m_1},\qquad
C=\sqrt{2/m_1}(A_0+C_0).
\label{eq:explicit_ES_constants}
\end{align}
Here $m_1,M_1$ are given in \eqref{eq:real_beam_final_norm_equivalence}, $M_L,\lambda_L$ in \eqref{eq:explicit_semigroup_constants}, and $A_0,A_1,K$ in \eqref{eq:explicit_normal_form_constants}.
\end{thm}
\endgroup

\begin{pf}
Fix $\omega\geq\omega_0$. By Proposition~\ref{prop:filtered_average_stability} and \eqref{eq:explicit_semigroup_constants}, the averaged generator in \eqref{eq:real_semilinear_ES_system} satisfies $\|e^{\mathcal Lt}\|\leq M_Le^{-\lambda_Lt}$. The perturbation \eqref{eq:real_zero_mean_perturbation} contains only the harmonics $0<|m|\leq4$. Appendix~\ref{app:beam_averaging} establishes
\begin{align}
\|(i\Omega-\mathcal L)^{-1}\mathcal B\|
\leq C_R|\Omega|^{-1/2},\qquad |\Omega|\geq1,
\label{eq:real_beam_resolvent_bound}
\end{align}
with $C_R$ from \eqref{eq:explicit_resolvent_constant}. Thus the bound applies at $\Omega=m\omega$ even when $|m|\omega$ equals a modal frequency of the PDE. The periodic change $Y=W+\mathcal R_\omega(\omega t,W)$ in \eqref{eq:real_beam_averaging_change} therefore yields \eqref{eq:real_beam_averaged_normal_form}. On $\|W\|\leq r_*$, the appendix gives, uniformly in the phase,
\begin{equation}
\begin{aligned}
\|\mathcal R_\omega\|&\leq A_0/\sqrt\omega,\\
\|D_W\mathcal R_\omega\|&\leq A_1/\sqrt\omega\leq\tfrac12,\\
\|\mathcal Q_\omega\|_{C^1}&\leq K/\sqrt\omega.
\end{aligned}
\label{eq:periodic_change_remainder_bounds}
\end{equation}
The change applies to mild solutions and has Lipschitz constant at most $3/2$. Its inverse is defined on $\|Y\|<r_*/2$ and has Lipschitz constant at most $2$.

On the space of continuous $\Pi$-periodic functions with the supremum norm, consider
\begin{align}
(\mathcal F W)(t)
=\int_0^\infty e^{\mathcal Ls}
\mathcal Q_\omega(\omega(t-s),W(t-s))\,ds.
\label{eq:real_beam_periodic_fixed_point}
\end{align}
For $\|W\|_\infty\leq r_*$, the semigroup estimate gives $\|\mathcal F W\|_\infty\leq C_0/\sqrt\omega$ and $\operatorname{Lip}(\mathcal F)\leq C_0/\sqrt\omega\leq1/2$. Since $C_0/\sqrt\omega\leq r_*/4$ by \eqref{eq:explicit_ES_constants}, $\mathcal F$ is a contraction on the closed ball of radius $r_*$. Its fixed point $W^\Pi$ is a periodic mild solution. Conversely, any periodic mild solution in this ball satisfies \eqref{eq:real_beam_periodic_fixed_point} by variation of constants from $t-N\Pi$ and passage to the limit $N\to\infty$.

Returning to $Y$ and using \eqref{eq:real_beam_final_norm_equivalence} gives a periodic solution $Z^\Pi$ with
\begin{equation}
\sup_t\mathcal N_{\rm EB}(Z^\Pi(t))
\leq \frac{A_0+C_0}{\sqrt{m_1\omega}}
\leq \frac{r_*}{4\sqrt{M_1}}.
\label{eq:periodic_error_neighborhood}
\end{equation}
This proves \eqref{eq:convergance_closedloop}. Any periodic solution in the neighborhood stated in the theorem satisfies $\|Y\|_\infty<r_*/2$. The periodic inverse change then gives a solution in $\|W\|_\infty\leq r_*$, where the fixed point is unique.

Now suppose $\mathcal N_{\rm EB}(Z(0)-Z^\Pi(0))<\delta$. The preceding estimate and \eqref{eq:explicit_ES_constants} give $\|Y(0)\|<r_*/4+\sqrt{M_1}\delta<r_*/2$. Hence the inverse change gives $\|W(0)-W^\Pi(0)\|\leq2\sqrt{M_1}\mathcal N_{\rm EB}(Z(0)-Z^\Pi(0))$. As long as $\|W(t)\|\leq r_*$, variation of constants and Gronwall's inequality yield
\begin{align}
\|W(t)-W^\Pi(t)\|
&\leq M_L e^{-(\lambda_L-M_LK/\sqrt\omega)t}\nonumber\\
&\qquad{}\times\|W(0)-W^\Pi(0)\|.
\label{eq:ES_periodic_attraction}
\end{align}

By \eqref{eq:explicit_ES_constants}, the decay rate is at least $\lambda_L/2$, and $2M_L\sqrt{M_1}\delta=r_*/4$. Since $\|W^\Pi\|_\infty\leq r_*/4$, it follows that $\|W(t)\|\leq r_*/2$ and $\|Y(t)\|\leq3r_*/4$. Using $|\vartheta_1|\leq d\|Y\|$ from Appendix~\ref{app:beam_averaging}, we obtain
\begin{equation}
|\Theta(t)-\Theta^*|\leq |a|+3dr_*/4<r_Q.
\label{eq:quadratic_region_continuation}
\end{equation}
These bounds allow continuation for all $t\geq0$. Applying the forward change and \eqref{eq:real_beam_final_norm_equivalence} to \eqref{eq:ES_periodic_attraction} gives \eqref{eq:uniform_physical_orbit_attraction} with $M_*=3M_L\sqrt{M_1/m_1}$ and $\lambda_*=\lambda_L/2$.

Finally, $\max\{|h(0)|,|h(1)|\}\leq\sqrt2\|h\|_{H^1}$ gives
\begin{equation}
\sup_t\max\{|\alpha^\Pi(t,0)|,|\alpha^\Pi(t,1)|\}
\leq C/\sqrt\omega
\label{eq:periodic_endpoint_bound}
\end{equation}
by the preceding bound on $\mathcal N_{\rm EB}(Z^\Pi)$ and the choice of $C$ in \eqref{eq:explicit_ES_constants}. By \eqref{eq:uniform_physical_orbit_attraction}, $\alpha-\alpha^\Pi$ converges to zero at both endpoints. From \eqref{eq:errorstate}, \eqref{eq:trajectorybc03}, and
\eqref{eq:perturbation1},
\begin{equation}
\begin{aligned}
\Theta(t)-\Theta^*&=\alpha(t,0)+a\sin(\omega t),\\
\theta_1(t)-\Theta^*&=\alpha(t,1)+S_1(t).
\end{aligned}
\label{eq:physical_endpoint_identities}
\end{equation}

Since $\sup_t|S_1(t)|=r_1(\omega)$ by \eqref{eq:perturbationS1}, these identities give \eqref{eq:convergence_Theta} and \eqref{eq:convergence_theta1}. The local quadratic map \eqref{extrachapter.eq:static_map} then gives \eqref{eq:convergence_y}. \hfill $\square$
\end{pf}

\begin{rem}\label{rem:beam_regularity_and_rate}
By \eqref{eq:perturbationS1}--\eqref{eq:perturbationS2}, the probing amplitudes in the position and moment channels are
\begin{align}
r_1(\omega)&=\tfrac{|a|}{2}|\cosh\sqrt\omega+\cos\sqrt\omega|,\nonumber\\
r_2(\omega)&=\tfrac{|a|\omega}{2}|\cosh\sqrt\omega-\cos\sqrt\omega|.
\label{eq:both_probing_amplitudes}
\end{align}
They satisfy $r_1\sim |a|e^{\sqrt\omega}/4$ and $r_2\sim |a|\omega e^{\sqrt\omega}/4$. Thus, improved error bounds entail increased probing effort.

\end{rem}

\section{Simulation Results}\label{section:simulation}

We simulate the system \eqref{eq:euler_bernoulli_pde}--\eqref{eq:euler_bernoulli_bc1} with the controllers \eqref{eq:control_law_filtered}--\eqref{eq:control_law_filtered_2}, using a finite-element discretization with $64$ cubic Hermite elements. Time integration uses the implicit trapezoidal rule with $\Delta t\approx5\times10^{-5}$, over a simulation horizon of $140$ normalized time units. For the quadratic map \eqref{extrachapter.eq:static_map}, we take $H=-1$, $\Theta^\ast=1.5$, and $y^\ast=2.4$. The probing amplitude and frequency are $a=0.2$ and $\omega=20$, respectively. The controller gains are $c=0.1$, $\overline c=2$, $\mu_p=0.02$, and $\mu_d=0.05$, satisfying \eqref{eq:filter_cutoff_condition} since $\overline c>\mu_p/\mu_d=0.4$. The initial conditions are $u(0,x)=0$, $u_t(0,x)=R_t(0,x)$, and $\eta(0)=F(0)=0$, consistently with \eqref{eq:ES_filter_initial_compatibility}. The initial velocity is represented through its cubic Hermite interpolant. Hence, the beam is initially undeformed but is not initially at rest.

\begin{figure*}[htpb]
    \centering
    \includegraphics[width=0.8\linewidth]{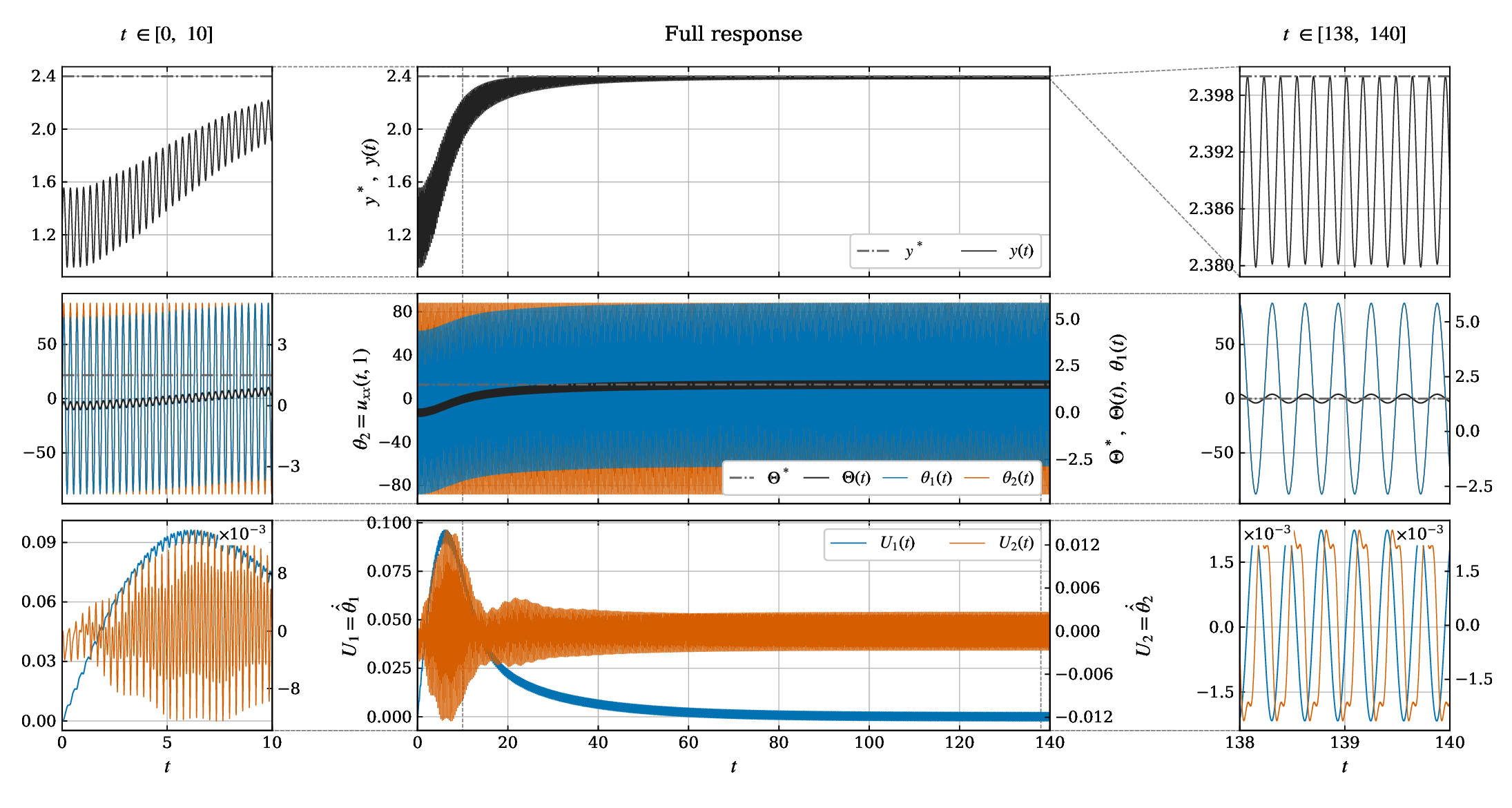}
    \caption{Extremum-seeking response: measured objective $y(t)$ and optimal value $y^*$ (top); map input $\Theta(t)$, optimal input $\Theta^*$, boundary displacement $\theta_1(t)$, and curvature $\theta_2(t)$ (middle); control signals $U_1(t)$ and $U_2(t)$ (bottom).}
    \label{fig:results}
\end{figure*}

Figure~\ref{fig:results} shows the measured objective, the map input, the boundary signals, and the nominal control signals. The center column presents the full responses, and the left and right columns give enlarged views of the intervals $t\in[0,10]$ and $t\in[138,140]$, respectively.
The top row shows the measured objective $y(t)$ approaching a neighborhood of its maximum $y^*=2.4$. For the prescribed motion $\Theta(t)=\Theta^*+0.2\sin(20t)$, the quadratic map gives $y(t)=2.4-0.02\sin^2(20t)$. This accounts for the approximate range $[2.38,2.40]$ and the two oscillations per probing period in the top-right panel.
In the middle row, the right axis shows the map input $\Theta(t)=u(t,0)$, the optimizer $\Theta^*$, and the actuated-end displacement $\theta_1(t)=u(t,1)$. The left axis shows the actuated-end curvature $\theta_2(t)=u_{xx}(t,1)$. Producing the prescribed amplitude $a=0.2$ at the map input requires position and moment probes of amplitudes $r_1(20)\approx4.354$ and $r_2(20)\approx88.031$, respectively, by \eqref{eq:perturbationS1}--\eqref{eq:perturbationS2}.
The persistent displacement oscillations in the middle-right panel are consistent with \eqref{eq:convergence_Theta}--\eqref{eq:convergence_theta1}. Under the hypotheses of Theorem~\ref{thm:ES_stability_revised}, the error relative to the prescribed probing motion converges to a periodic trajectory of size $\mathcal O(\omega^{-1/2})$.
The bottom row shows the nominal control signals $U_1=\dot{\hat{\theta}}_1$ (left axis) and $U_2=\dot{\hat{\theta}}_2$ (right axis). Both exhibit oscillations on the $10^{-3}$ scale in the bottom-right panel.

\begin{figure*}[htpb]
    \centering
    \includegraphics[width=0.8\linewidth]{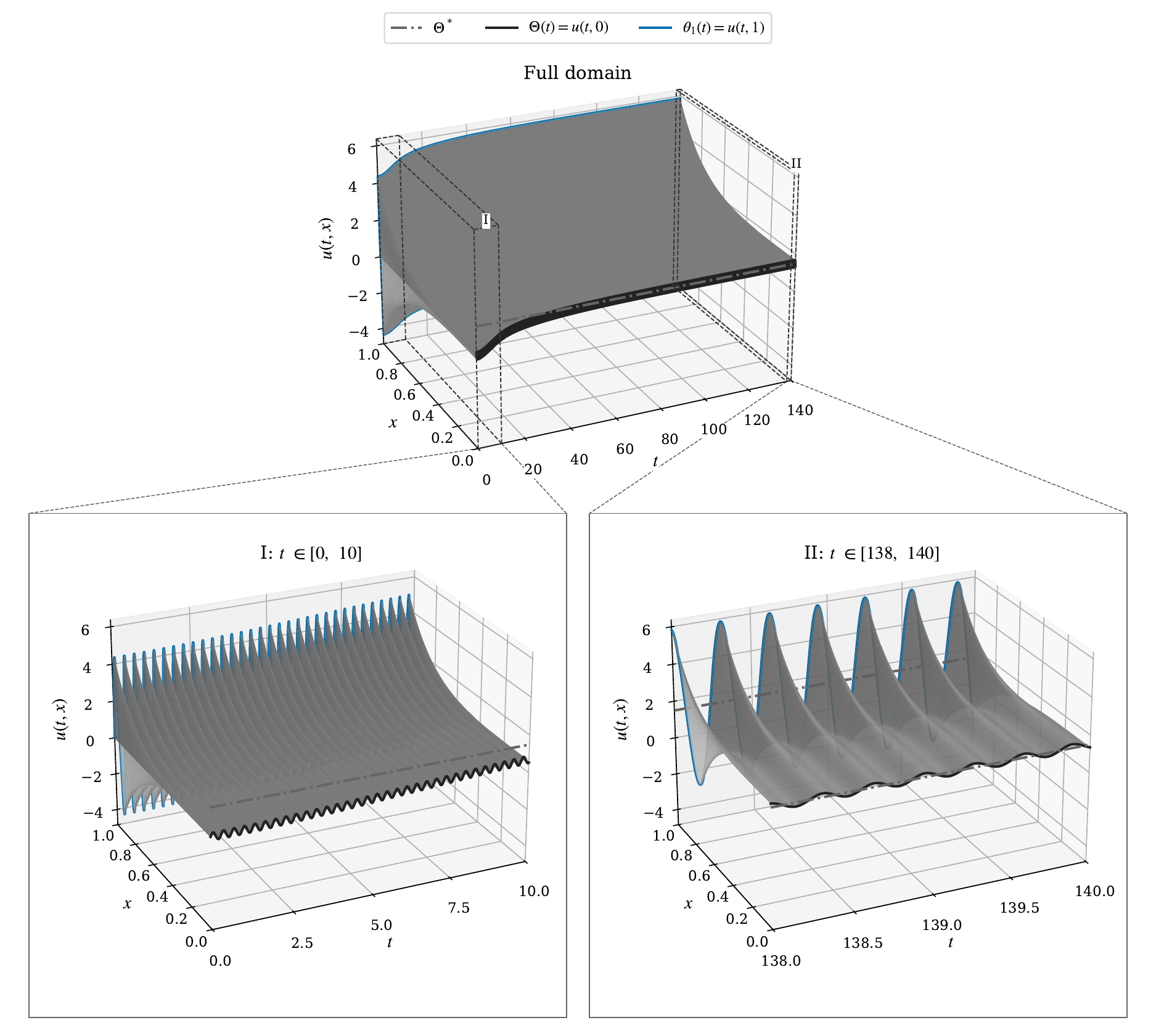}
    \caption{Displacement $u(t,x)$ over the full simulation and the two indicated time intervals. The dashed line marks $\Theta^\ast$ at $x=0$.}
    \label{fig:results_3D}
\end{figure*}

Figure~\ref{fig:results_3D} presents the spatiotemporal evolution of the beam displacement $u(t,x)$ over the normalized spatial domain $x\in[0,1]$. As can be seen in the detailed zoom in Figure~\ref{fig:results_3D}, the sliding-end displacement oscillates about $\Theta^\ast$ with amplitude close to $0.2$, while the controlled-end displacement oscillates with amplitude close to $4.35$.

\section{Conclusions} \label{section:conclusion}
This paper developed a boundary ES controller for an unknown locally quadratic map of the sliding-end displacement of an Euler--Bernoulli beam. The feedback uses the known beam model, distributed acceleration $u_{tt}$, and curvature rate $u_{txx}$. Only the map-dependent terms pass through the low-pass filter. The averaged closed loop is exponentially stable for $\overline c>\mu_p/\mu_d$. For a fixed probing amplitude within the quadratic region and sufficiently high probing frequency, the full closed loop admits a locally exponentially attracting periodic error trajectory. The resulting bounds distinguish the accuracy at the map input from the boundary probing amplitudes, clarifying the trade-off between optimization accuracy and probing effort.

Future research directions include extending the design to other beam boundary conditions and actuator configurations, developing observer-based implementations to reduce the measurement requirements, and incorporating actuator constraints to account for the boundary probing demands.



\appendix
\renewcommand*{\theHsection}{appendix.\Alph{section}}

\section{Physical Reconstruction and Well-Posedness}\label{app:beam_coordinates}

We derive the real coordinates for the exact error system, prove \eqref{eq:real_beam_final_norm_equivalence}, and identify the generator of the averaged closed loop.

\subsection{Real coordinates for the closed-loop error system}

Represent each complex number by its real part and negative imaginary part. In these coordinates, the Volterra operator in \eqref{eq:backsteppingtransformation2} is
\begin{align}
(\mathcal T_c g)(x)&=g(x)-\int_0^x K_c(x,y)g(y)\,dy,\nonumber\\
K_c(x,y)&=\begin{pmatrix}
\kappa_r(x,y)&\kappa_i(x,y)\\
-\kappa_i(x,y)&\kappa_r(x,y)
\end{pmatrix}.
\label{eq:real_Volterra_operator}
\end{align}
Set $\mathbf b=(b_1,b_2)=(-\alpha_{xxxx},\beta_{xx})$ and $(r_c,s_c)=\mathcal T_c\mathbf b$. For compatible classical solutions, $\beta_t=-\alpha_{xxxx}$, so $\mathbf b$ represents $v$ in \eqref{eq:transformation_schrodinger}. The kernel equations \eqref{eq:kernelpdefinal}--\eqref{eq:kernelbc2final} and the feedback \eqref{eq:split_filtered_complex_control} give
\begin{align}
(r_c)_t&=-(s_c)_{xx}-cr_c,\nonumber\\
(s_c)_t&=(r_c)_{xx}-cs_c,\nonumber\\
(r_c)_x(t,0)&=(s_c)_x(t,0)=0,\nonumber\\
(r_c(t,1),s_c(t,1))&=f(t).
\label{eq:real_beam_field_dynamics}
\end{align}

Define the four scalar coordinates entering the filter by
\begin{align}
\eta_b(t,x)&=\eta(t)-\int_0^x(x-y)\beta(t,y)\,dy,\nonumber\\
\mathbf z_1=(z_1,z_3)&=[\mathcal T_c(\alpha,\eta_b)](1),\nonumber\\
\mathbf z_2=(z_2,z_4)&=[\mathcal T_c(\beta,\alpha_{xx})](1).
\label{eq:real_beam_predictors}
\end{align}
By \eqref{eq:eta_extension} and \eqref{eq:erroTheta}, $\dot\eta=\alpha_{xx}(t,0)$. Integrating $\beta_t=-\alpha_{xxxx}$ twice and using $\alpha_{xxx}(t,0)=0$ gives $(\eta_b)_t=\alpha_{xx}$. Hence \eqref{eq:real_beam_predictors} and the boundary condition at $x=1$ give
\begin{align}
\dot{\mathbf z}_1=\mathbf z_2,\qquad \dot{\mathbf z}_2=f.
\label{eq:real_predictor_dynamics}
\end{align}
With $k_p=-H\mu_p$ and $k_d=-H\mu_d$, the kernels \eqref{eq:solutionkernelode_intermediate}--\eqref{eq:solutionkernelpde_intermediate} give $z=v-H[\mu_p(\alpha-i\eta_b)+\mu_d(\beta-i\alpha_{xx})]$ in \eqref{eq:backsteppingtransformation1}. Applying \eqref{eq:backsteppingtransformation2} yields
\begin{align}
(r_c,s_c)&=(\operatorname{Re}w,-\operatorname{Im}w)\nonumber\\
&\quad+H\mathcal T_c\bigl[\mu_p(\alpha,\eta_b)
                         +\mu_d(\beta,\alpha_{xx})\bigr].
\label{eq:full_target_energy_coordinate_bridge}
\end{align}
At $x=1$, $r_c-is_c=F$. Subtracting the deterministic feedback \eqref{eq:complex_control_law_mean} from the applied input \eqref{eq:split_filtered_complex_control} gives $w(t,1)=F-H\Phi$ before averaging. Thus
\begin{equation}
\Phi=\mu_p(z_1-iz_3)+\mu_d(z_2-iz_4).
\label{eq:exact_predictor_kernel_identity}
\end{equation}
With $h_1^{\rm R}=(\operatorname{Re}\bar h_1,-\operatorname{Im}\bar h_1)$, the exact filter \eqref{eq:exact_complex_closed_loop_filter} and the averaged filter \eqref{eq:averaged_complex_filter} become
\begin{align}
\dot f&=-\overline c f+\overline c\hat H(\mu_p\mathbf z_1+\mu_d\mathbf z_2)
\nonumber\\
&\quad+\overline c h_1^{\rm R}(G-\hat H\vartheta_1),
\label{eq:exact_predictor_ES_signal}\\
\dot f_{\rm av}&=-\overline c f_{\rm av}
+\overline c H(\mu_p\mathbf z_{1,{\rm av}}+\mu_d\mathbf z_{2,{\rm av}}).
\label{eq:real_predictor_feedback_identity}
\end{align}

To obtain homogeneous boundary conditions, let $(\mathcal I g)(x)=-\int_x^1\int_0^y g(\xi)\,d\xi\,dy$, so that $(\mathcal I g)''=g$, $(\mathcal I g)'(0)=(\mathcal I g)(1)=0$. Set $k=\sqrt{c/2}$, $A(x)=\cosh(kx)\cos(kx)$ and $B(x)=\sinh(kx)\sin(kx)$. The profiles
\begin{align}
a_c(x)&=\frac{A(x)A(1)+B(x)B(1)}{A(1)^2+B(1)^2},\nonumber\\
b_c(x)&=\frac{B(x)A(1)-A(x)B(1)}{A(1)^2+B(1)^2}
\label{eq:real_lifting_explicit}
\end{align}
satisfy $a_c''=-cb_c$, $b_c''=ca_c$, $a_c'(0)=b_c'(0)=0$ and $(a_c(1),b_c(1))=(1,0)$. The change of variables
\begin{align}
\chi&=\mathcal I(s_c+b_cf_1-a_cf_2),\nonumber\\
\upsilon&=r_c-a_cf_1-b_cf_2
\label{eq:real_beam_lifting}
\end{align}
therefore transforms \eqref{eq:real_beam_field_dynamics} into
\begin{align}
\chi_t&=\upsilon-c\chi+\mathcal I(b_c\dot f_1-a_c\dot f_2),\nonumber\\
\upsilon_t&=-\chi_{xxxx}-c\upsilon-a_c\dot f_1-b_c\dot f_2,
\label{eq:lifted_real_beam_dynamics}\\
\chi_x(t,0)=\chi_{xxx}(t,0)&=\chi(t,1)=\chi_{xx}(t,1)=0,\nonumber\\
\upsilon_x(t,0)&=\upsilon(t,1)=0.
\label{eq:lifted_real_beam_boundaries}
\end{align}

\subsection{Reconstruction and norm equivalence}

Let $e_1=(\vartheta_1,\eta)$, $e_2=(p,\vartheta_2)$, and define $(\mathcal V_j g)(x)=\int_0^x(x-y)^{j-1}g(y)\,dy/(j-1)!$ for $j\geq1$. Integration from the sliding end gives
\begin{align}
\beta(t,x)&=p(t)+(\mathcal V_2 b_2)(t,x),\nonumber\\
\alpha(t,x)&=\vartheta_1(t)+\tfrac{x^2}{2}\vartheta_2(t)
             -(\mathcal V_4 b_1)(t,x).
\label{eq:real_beam_reconstruction}
\end{align}
The trace bound $|h(0)|^2\leq2\|h\|_{H^1}^2$ \citep{AdamsFournier2003}, applied to $\alpha,\alpha_{xx},\beta$, yields the upper bound in \eqref{eq:physical_reconstruction_norm_bounds}. For the lower bound, differentiate \eqref{eq:real_beam_reconstruction} and use $\|\mathcal V_j\|_{L^2\to L^2}\leq1/j!$. Together these estimates give
\begin{align}
\tfrac13\mathcal N_{\rm EB}(Z)^2
&\leq |e_1|^2+|e_2|^2+\|\mathbf b\|_{L^2}^2+|f|^2\nonumber\\
&\leq2\mathcal N_{\rm EB}(Z)^2.
\label{eq:physical_reconstruction_norm_bounds}
\end{align}

Given $(\xi,\chi,\upsilon)$, recover $f$ from \eqref{eq:real_filter_matrix}, then $r_c=\upsilon+a_cf_1+b_cf_2$ and $s_c=\chi_{xx}-b_cf_1+a_cf_2$ from \eqref{eq:real_beam_lifting}. Thus $\mathbf b=\mathcal T_c^{-1}(r_c,s_c)$. Substitution of \eqref{eq:real_beam_reconstruction} into \eqref{eq:real_beam_predictors} gives
\begin{align}
\mathbf z_2&=D_0e_2+\mathcal J_2\mathbf b,\nonumber\\
\mathbf z_1&=D_0e_1+D_2e_2+\mathcal J_1\mathbf b,
\label{eq:real_predictor_triangular_inverse}
\end{align}
where $J=\left(\begin{smallmatrix}0&1\\-1&0\end{smallmatrix}\right)$ and
\begin{align}
D_0&=\begin{pmatrix}A(1)&B(1)\\-B(1)&A(1)\end{pmatrix},\qquad
D_2=\tfrac12[\mathcal T_c(x^2J)](1),\nonumber\\
\mathcal J_1\mathbf b&=-[\mathcal T_c\mathcal V_4\mathbf b](1),\qquad
\mathcal J_2\mathbf b=[\mathcal T_c J\mathcal V_2\mathbf b](1).
\label{eq:predictor_inverse_coefficients}
\end{align} 
Here $\mathcal V_j$ acts componentwise and $\mathcal T_c$ acts on matrix columns. The kernel equations give $P_c''=icP_c$, $P_c(0)=1$, $P_c'(0)=0$ for $P_c(x)=1-\int_0^x\kappa(x,y)\,dy$. Thus $P_c=A+iB$, which gives $D_0$. Since $A(x)^2+B(x)^2=[\cosh(2kx)+\cos(2kx)]/2$ increases from $1$, $\|D_0^{-1}\|\leq1$ and $a_c^2+b_c^2\leq1$. Equation~\eqref{eq:real_predictor_triangular_inverse} therefore determines $e_2$ and then $e_1$, completing the reconstruction.

Set $K_\kappa=(c/2)e^{c/4}$, which bounds $|\kappa|$ on $0\leq y\leq x\leq1$ by the power series of $I_1(z)/z$. Let $L=1+K_\kappa$ and $E=e^{K_\kappa}$. The Volterra series gives $\|\mathcal T_c\|\leq L$ and $\|\mathcal T_c^{-1}\|\leq E$ on $L^2$. The endpoint estimate $|(\mathcal T_c g)(1)|\leq L\|g\|_\infty$ gives $\|D_0\|,\|D_2\|,\|\mathcal J_1\|,\|\mathcal J_2\|\leq L$ by \eqref{eq:predictor_inverse_coefficients}. Writing $b=\|\mathbf b\|_{L^2}$ and $s=\|(\upsilon,\chi_{xx})\|_{L^2}$,
\eqref{eq:real_beam_lifting} gives $s\leq Lb+|f|$ and $b\leq E(s+|f|)$. Solving \eqref{eq:real_predictor_triangular_inverse} successively gives $|e_2|\leq|\mathbf z_2|+Lb$ and $|e_1|\leq|\mathbf z_1|+L|\mathbf z_2|+(L+L^2)b$. Squaring and summing the forward and inverse bounds by Cauchy--Schwarz gives $8L^2$ and $16E^2L^4$ for the squared norms of the map $(e_1,e_2,\mathbf b,f)\mapsto(\xi,\chi_{xx},\upsilon)$ and its inverse, respectively. With \eqref{eq:physical_reconstruction_norm_bounds},
this proves \eqref{eq:real_beam_final_norm_equivalence} with
\begin{align}
m_1=\frac{e^{-2K_\kappa}}{48(1+K_\kappa)^4},\qquad
M_1=16(1+K_\kappa)^2.
\label{eq:explicit_coordinate_constants}
\end{align}

\subsection{Generator of the averaged closed loop}

Let $\mathscr B=H_b^2\times L^2$, where
\begin{align}
H_b^2&=\{h\in H^2:h'(0)=h(1)=0\},\nonumber\\
\|(h,q)\|_{\mathscr B}^2&=\|h''\|_{L^2}^2+\|q\|_{L^2}^2.
\label{eq:real_beam_phase_space}
\end{align}
Since $\mathcal I:L^2\to H_b^2$ inverts $h\mapsto h''$, the reconstruction defines an isomorphism $\mathcal P:Z\mapsto(\xi,\chi,\upsilon)$ from $\mathcal X_{\rm EB}$ onto $\mathscr H_b=\mathbb R^6\times\mathscr B$. The homogeneous PDE in \eqref{eq:lifted_real_beam_dynamics} has generator
\begin{align}
\mathcal A_b(h,q)&=(q-ch,-h^{(4)}-cq),\nonumber\\
D(\mathcal A_b)&=\{(h,q)\in(H^4\cap H_b^2)\times H_b^2:\nonumber\\
&\hspace{25mm}h'''(0)=h''(1)=0\}.
\label{eq:real_fourth_order_generator}
\end{align}
In the basis $\sqrt2\cos(\nu_nx)$, $\nu_n=(n+\tfrac12)\pi$, $n\geq0$, each pair $(\nu_n^2h_n,q_n)$ evolves by a rotation scaled by $e^{-ct}$. Thus $\mathcal A_b$ generates a $C_0$-semigroup of norm $e^{-ct}$. The input operator
\begin{align}
\mathcal B_b g
=\bigl(\mathcal I(b_cg_1-a_cg_2),-a_cg_1-b_cg_2\bigr)
\label{eq:real_beam_input_operator}
\end{align}
satisfies $\|\mathcal B_b g\|_{\mathscr B}^2 =|g|^2\int_0^1(a_c^2+b_c^2)\,dx=D_c^2|g|^2$. For the averaged filter, $\dot f=L_f\xi$ as in the proof of Proposition~\ref{prop:filtered_average_stability}. Omitting the subscript $\mathrm{av}$ gives
\begin{align}
\dot\xi&=M_f\xi,\nonumber\\
\partial_t(\chi,\upsilon)&=\mathcal A_b(\chi,\upsilon)+\mathcal B_bL_f\xi.
\label{eq:real_lifted_average_system}
\end{align}
Since $\mathcal B_bL_f$ is bounded, the bounded perturbation theorem \citep[Theorem~III.1.3]{EngelNagel2000} gives a $C_0$-semigroup with generator $\mathcal L$ on $\mathscr H_b$ and domain $D(\mathcal L)=\mathbb R^6\times D(\mathcal A_b)$. By similarity,
\begin{align}
\mathcal A_{\rm EB}&=\mathcal P^{-1}\mathcal L\mathcal P,\qquad
D(\mathcal A_{\rm EB})=\mathcal P^{-1}D(\mathcal L).
\label{eq:physical_generator_domain}
\end{align}
In the original coordinates, this domain consists of $Z\in\mathcal X_{\rm EB}$ with $\alpha\in H^6$, $\beta\in H^4$, $\alpha^{(5)}(0)=\beta'''(0)=0$, and $(\mathcal T_c\mathbf b)(1)=f$. $\mathcal T_c$ and its inverse preserve $H^2$, and $K_c(0,0)=0$ gives $(\mathcal T_c\mathbf b)_x(0)=\mathbf b_x(0)$. The condition at $x=1$ follows from \eqref{eq:real_beam_lifting} and agrees with \eqref{eq:ES_filter_initial_compatibility} at $t=0$. Similarity gives well-posedness on $\mathcal X_{\rm EB}$, and density of $D(\mathcal A_{\rm EB})$ extends \eqref{eq:physical_EB_exponential_stability} to every mild solution.

\section{Averaging Estimates}\label{app:beam_averaging}

In the coordinates $Y=(\xi,\chi,\upsilon)$ of Appendix~\ref{app:beam_coordinates}, the operators in \eqref{eq:real_semilinear_ES_system} are
\begin{align}
\mathcal LY&=(M_f\xi,\mathcal A_b(\chi,\upsilon)+\mathcal B_bL_f\xi),
\nonumber\\
\mathcal Bg&=(E_fg,\mathcal B_bg),\qquad E_fg=(0,0,g_1,0,0,g_2),
\label{eq:real_augmented_beam_operators}
\end{align}
with $D(\mathcal L)=\mathbb R^6\times D(\mathcal A_b)$. Subtracting \eqref{eq:real_predictor_feedback_identity} from \eqref{eq:exact_predictor_ES_signal} at the same state gives
\begin{align}
\Delta(\tau,\mathbf r(Y))
&=\overline c(\hat H-H)(\mu_p\mathbf z_1+\mu_d\mathbf z_2)\nonumber\\
&\quad+\overline c h_1^{\rm R}(G-\hat H\vartheta_1).
\label{eq:real_zero_mean_perturbation}
\end{align}
For fixed $Y$ with $|\vartheta_1|+|a|<r_Q$, the quadratic map and the demodulation formulas \eqref{eq:gradient_estimate} and \eqref{eq:hessian} give harmonics up to order four in $\tau=\omega t$. The constant term in \eqref{eq:real_zero_mean_perturbation} vanishes by \eqref{eq:hessian_estimates}. Smooth state dependence gives local well-posedness \citep[Chapter~6]{Pazy1983}. By
\eqref{eq:real_predictor_dynamics} and $\dot\vartheta_1=p$,
\begin{equation}
\mathbf r(\mathcal LY)=\mathbf r_L(Y)=(z_2,f_1,z_4,f_2,p),
\quad Y\in D(\mathcal L).
\label{eq:ES_error_vector_derivative}
\end{equation}
Applying the endpoint trace bound to $\vartheta_1=\alpha(t,0)$ and $p=\beta(t,0)$ and using \eqref{eq:real_beam_final_norm_equivalence} gives $\max\{|\mathbf r(Y)|,|\mathbf r_L(Y)|\}\leq d\|Y\|$ on $\mathscr H_b$, where $d=\sqrt{1+2/m_1}$.

\smallskip\noindent\textit{Resolvent estimate.}
Let $s=|\Omega|\geq1$ and $R_b=(i\Omega I-\mathcal A_b)^{-1}$. In the cosine basis of Appendix~\ref{app:beam_coordinates}, integration by parts in \eqref{eq:real_lifting_explicit} gives the coefficients of $a_c,b_c$ as
\begin{equation}
a_n+ib_n=\frac{\sqrt2(-1)^n\nu_n}{\nu_n^2+ic},
\qquad \nu_n=(n+\tfrac12)\pi.
\label{eq:exact_lifting_coefficients}
\end{equation}
Since $a_n^2+b_n^2\leq2/\nu_n^2$ and the modal matrices are normal with eigenvalues $-c\pm i\nu_n^2$, \eqref{eq:real_beam_input_operator} gives
\begin{align}
\|R_b\mathcal B_b\|^2
&\leq\sum_{n\geq0}\frac{2/\nu_n^2}{c^2+(s-\nu_n^2)^2}
\leq\frac{C_b^2}{s},\nonumber\\
C_b^2&=4+4/c^2+1/\pi^2.
\label{eq:real_beam_modal_resolvent}
\end{align}
The terms with $\nu_n^2\leq s/2$ sum to at most $4/s^2$, since $\sum_n\nu_n^{-2}=1/2$. Elsewhere $2/\nu_n^2\leq4/s$. Counting from the nearest modal frequency and using the minimum spacing $d_0=2\pi^2$ gives
\begin{equation}
\begin{aligned}
\sum_{n\geq0}\frac{1}{c^2+(s-\nu_n^2)^2}
&\leq \frac1{c^2}+\frac{2}{d_0^2}\sum_{k\geq1}(k-\tfrac12)^{-2}\\
&=\frac1{c^2}+\frac1{4\pi^2}.
\end{aligned}
\label{eq:modal_reciprocal_sum}
\end{equation}
For $R_M=(i\Omega I-M_f)^{-1}$, the identities $J_fE_f=I_2$ and $M_fR_M=i\Omega R_M-I$ reduce the full resolvent to
\begin{align}
(i\Omega I-\mathcal L)^{-1}\mathcal Bg
=\bigl(R_ME_fg,\,i\Omega R_b\mathcal B_bJ_fR_ME_fg\bigr).
\label{eq:augmented_resolvent_blocks}
\end{align}
Let $q_{\min},q_{\max}$ denote the extreme eigenvalues of the Lyapunov matrix $Q$ in Proposition~\ref{prop:filtered_average_stability}. For $\Omega_c=\max\{1,2\|M_f\|\}$, its Lyapunov bound for $s\leq\Omega_c$ and the Neumann series for $s\geq\Omega_c$ give $s\|R_M\|\leq C_M$. Thus \eqref{eq:real_beam_modal_resolvent} and \eqref{eq:augmented_resolvent_blocks} prove \eqref{eq:real_beam_resolvent_bound} with
\begin{align}
C_M&=\max\{2,\,2\Omega_c q_{\max}\sqrt{q_{\max}/q_{\min}}\},
\nonumber\\
C_R&=C_M\sqrt{1+C_b^2}.
\label{eq:explicit_resolvent_constant}
\end{align}

\smallskip\noindent\textit{Periodic change of coordinates.}
Let $\Delta_m$ be the Fourier coefficients of $\Delta$ and $S_m=(im\omega I-\mathcal L)^{-1}\mathcal B$, with sums over $0<|m|\leq4$. Conjugate harmonics give the real periodic change
\begin{align}
Y&=W+\mathcal R_\omega(\omega t,W),\nonumber\\
\mathcal R_\omega(\tau,W)&=\sum_m e^{im\tau}S_m\Delta_m(\mathbf r(W)).
\label{eq:real_beam_averaging_change}
\end{align}
At fixed $W$, $\omega\partial_\tau\mathcal R_\omega
=\mathcal L\mathcal R_\omega+\mathcal B\Delta(\tau,\mathbf r(W))$.
Substitution in \eqref{eq:real_semilinear_ES_system} yields
\begin{align}
\dot W&=\mathcal LW+\mathcal Q_\omega(\omega t,W),\nonumber\\
\mathcal Q_\omega&=(I+D_W\mathcal R_\omega)^{-1}N_\omega,\nonumber\\
N_\omega&=\mathcal B[\Delta(\tau,\mathbf r(W+\mathcal R_\omega))
                     -\Delta(\tau,\mathbf r(W))]\nonumber\\
&\quad-D_W\mathcal R_\omega\,\mathcal LW.
\label{eq:real_beam_averaged_normal_form}
\end{align}

By \eqref{eq:ES_error_vector_derivative}, the last term is defined throughout the local state space by the identity
\begin{equation}
D_W\mathcal R_\omega\,\mathcal LW
=\sum_m e^{im\tau}S_m D\Delta_m(\mathbf r(W))\mathbf r_L(W),
\label{eq:normal_form_unbounded_term}
\end{equation}
initially valid on $D(\mathcal L)$.

Set $\varepsilon=\omega^{-1/2}$, $r_*=(r_Q-|a|)/(4d)$, and $B_0=\|\mathcal B\|=\sqrt{1+D_c^2}$. For $j=0,1,2$, define
\begin{align}
\Gamma_j&=\sum_m\max_{|u|\leq2dr_*}\|D^j\Delta_m(u)\|,
\qquad A_j=C_Rd^j\Gamma_j,\nonumber\\
F_0&=d\Gamma_1(B_0A_0+C_Rr_*),\nonumber\\
F_1&=d^2\Gamma_2(B_0A_0+C_Rr_*)+d\Gamma_1(B_0A_1+C_R),
\nonumber\\
K&=2(F_0+F_1)+4A_2F_0.
\label{eq:explicit_normal_form_constants}
\end{align}
On $\|W\|\leq r_*$, \eqref{eq:real_beam_resolvent_bound} gives $\|D_W^j\mathcal R_\omega\|\leq A_j\varepsilon$ uniformly in phase. For $\omega\geq\omega_0$, \eqref{eq:explicit_ES_constants} gives $A_0\varepsilon\leq r_*/4$ and $A_1\varepsilon\leq1/2$. Both arguments of $\Delta$ then lie in $|u|\leq2dr_*$. The mean value inequality and \eqref{eq:normal_form_unbounded_term} give $\|N_\omega\|\leq F_0\varepsilon$ and $\|D_WN_\omega\|\leq F_1\varepsilon$. The inverse factor in \eqref{eq:real_beam_averaged_normal_form} has norm at most $2$ and its state derivative has norm at most $4A_2\varepsilon$. Hence $\sup_\tau\|\mathcal Q_\omega(\tau,\cdot)\|_{C^1}\leq K\varepsilon$.

For $\|Y\|<r_*/2$, $W\mapsto Y-\mathcal R_\omega(\tau,W)$ is a contraction on the closed ball $\|W\|\leq r_*$. The forward and inverse changes have Lipschitz constants at most $3/2$ and $2$, respectively. Finally, $\mathcal LS_m=im\omega S_m-\mathcal B$ gives graph-norm regularity of the correction for fixed $\omega$. The change thus preserves $D(\mathcal L)$, and density and continuous dependence extend \eqref{eq:real_beam_averaged_normal_form} to mild solutions.

\textbf{Declaration of generative AI use.}
The authors used ChatGPT to assist with drafting, proof development, and simulation code, and Claude to review the theorem statements and proofs. The authors reviewed and edited the material and take full responsibility for the final manuscript.

\IfFileExists{references.bbl}{\input{references.bbl}}{%
\IfFileExists{bibli.bib}{%
  \bibliographystyle{elsarticle-harv}%
  \bibliography{bibli}%
}{%
  \PackageWarningNoLine{automatica-draft}{Upload bibli.bib to resolve the citations}%
}%
}

\raggedbottom

\end{document}